\documentclass[11pt, reqno]{amsart}

\usepackage{graphicx}

\usepackage{mathrsfs}
\usepackage{amssymb}
\usepackage{amsfonts, mathtools}
\usepackage{amssymb,amsmath,amsthm,tikz-cd, tikz}

\usepackage{hyperref}
\usepackage{palatino}
\usepackage{wrapfig}
\usepackage{subcaption}

\usetikzlibrary{patterns, matrix,shapes,arrows,calc,topaths,intersections,positioning,decorations.pathreplacing,decorations.pathmorphing,fit}

\newif\ifpaper
\papertrue

\usepackage{color}

\def\name{Peng Zhou}
\hypersetup{
  colorlinks = true,
  urlcolor = black,
  pdfauthor = {\name},
  pdfkeywords = {mathematics},
  pdftitle = {\name: },
  pdfsubject = {Mathematics},
  pdfpagemode = UseNone
}

\newcommand*{\hrlen}{5}
\newcommand*{\hramp}{3}

\tikzset{
asdstyle/.style={blue,thick},
righthairs/.style={postaction={decorate,draw,decoration={border,amplitude=\hramp,segment length=\hrlen,angle=-90,pre=moveto,pre length=\hrlen/2}}},
lefthairs/.style={postaction={decorate,draw,decoration={border,amplitude=\hramp,segment length=\hrlen,angle=90,pre=moveto,pre length=\hrlen/2}}},
righthairsnogap/.style={postaction={decorate,draw,decoration={border,amplitude=\hramp,segment length=\hrlen,angle=-90}}},
lefthairsnogap/.style={postaction={decorate,draw,decoration={border,amplitude=\hramp,segment length=\hrlen,angle=90}}},
graphstyle/.style={thick},
arrowstyle/.style={thick,decorate,decoration={snake,amplitude=1.7,segment length=10pt,post length=.5mm,pre length=0}},
genmapstyle/.style={thick,-stealth'},
arrhdstyle/.style={thick},
exceptarcstyle/.style={red, ultra thick},
dualquiverstyle/.style={thick,->},
patstyle/.style={pattern color = gray, pattern = north east lines, opacity=0.3}
}

\title {Lagrangian Skeleta of Hypersurfaces in $(\C^*)^n$}
\author{\name}
\thanks{This work is supported by an
IHES Simons Postdoctoral Fellowship as part of the Simons
Collaboration on HMS. }
\address{Institut des Hautes Études Scientifiques. Le Bois-Marie, 35 route de Chartres, 91440 Bures-sur-Yvette France}
\email{pengzhou@ihes.fr}
\date{\today}

\usepackage{color}

\renewcommand{\ss}{\subsection}
\newcommand{\z}{\text}

\newcommand{\Kahler}{K\"ahler }

\DeclareMathOperator{\cone}{cone}
\DeclareMathOperator{\conv}{conv}

\DeclareMathOperator{\Crit}{Crit}

\DeclareMathOperator{\dist}{dist}

\DeclareMathOperator{\Hess}{Hess}

\DeclareMathOperator{\Int}{Int}
\DeclareMathOperator{\id}{id}

\renewcommand{\Im}{\text{Im}}
\DeclareMathOperator{\Log}{Log}

\renewcommand{\Re}{\text{Re}}

\DeclareMathOperator{\spann}{span}

\newcommand{\pa}{\partial}

\newcommand{\la}{\langle}
\newcommand{\ra}{\rangle}
\newcommand{\ot}{\otimes}

\newcommand{\dbar}{\bar{\partial}}
\renewcommand{\d}{\partial}

\newcommand{\half}{\frac{1}{2}}

\newcommand{\wt}{\widetilde}
\newcommand{\wh}{\widehat}
\newcommand{\h}{\hat}
\newcommand{\wb}{\overline}

\newcommand{\acal}{\mathcal{A}}

\newcommand{\ccal}{\mathcal{C}}

\newcommand{\hcal}{\mathcal{H}}

\newcommand{\kcal}{\mathcal{K}}

\newcommand{\pcal}{\mathcal{P}}

\newcommand{\scal}{\mathcal{S}}
\newcommand{\tcal}{\mathcal{T}}

\newcommand{\wcal}{\mathcal{W}}

\newcommand{\C}{\mathbb{C}}

\newcommand{\R}{\mathbb{R}}

\newcommand{\Z}{\mathbb{Z}}

\newcommand{\isoto}{\xrightarrow{\sim}}
\newcommand{\lra}{\leftrightarrow}

\newcommand{\RA}{\Rightarrow}
\newcommand{\RM}{\backslash}

\newcommand{\into}{\hookrightarrow}

\newcommand{\onto}{\twoheadrightarrow}

\newcommand{\bea}{\begin{eqnarray*} }
\newcommand{\eea}{\end{eqnarray*} }
\newcommand{\be}{\begin{equation} }
\newcommand{\ee}{\end{equation} }
\newcommand{\bp}{\begin{proposition}}
\newcommand{\ep}{\end{proposition}}
\newcommand{\bt}{\begin{mmt}}
\newcommand{\et}{\end{mmt}}
\newcommand{\btu}{\begin{theou}}
\newcommand{\etu}{\end{theou}}
\newcommand{\bpf}{\begin{proof}}
\newcommand{\epf}{\end{proof}}
\newcommand{\bl}{\begin{lemma}}
\newcommand{\el}{\end{lemma}}
\newcommand{\bc}{\begin{corollary}}
\newcommand{\ec}{\end{corollary}}
\newcommand{\bd}{\begin{definition}}
\newcommand{\ed}{\end{definition}}
\newcommand{\bex}{\begin{example}}
\newcommand{\eex}{\end{example}}
\newcommand{\bA}{\left(\begin{array}}
\newcommand{\eA}{\end{array}\right)}
\newcommand{\bma}{\begin{bmatrix}}
\newcommand{\ema}{\end{bmatrix}}
\newcommand{\bcd}{\begin{tikzcd}}
\newcommand{\ecd}{\end{tikzcd}}
\newcommand{\bcs}{\begin{cases}}
\newcommand{\ecs}{\end{cases}}
\newcommand{\bee}{\begin{eqnarray} }
\newcommand{\eee}{\end{eqnarray} }
\newcommand{\brem}{\begin{remark}}
\newcommand{\erem}{\end{remark}}
\newcommand{\bnum}{\begin{enumerate}}
\newcommand{\enum}{\end{enumerate}}

\newtheorem{mmt}{Theorem}
\newtheorem*{theou}{Theorem}
\newtheorem*{maintheo}{Main Theorem}
\newtheorem{theo}{Theorem}[section]

\newtheorem{lemma}[theo]{Lemma}
\newtheorem{corollary}[theo]{Corollary}
\newtheorem{proposition}[theo]{Proposition}

\theoremstyle{definition}
\newtheorem{definition}[theo]{Definition}
\newtheorem{example}[theo]{Example}
\newtheorem{remark}[theo]{Remark}

\theoremstyle{plain}
\numberwithin{equation}{section}

\newcommand{\CS}{{\C^*}}
\newcommand{\CSN}{(\C^*)^n}

\begin{document}

\begin{abstract}
Let $W(z_1, \cdots, z_n): (\C^*)^n \to \C$ be a Laurent polynomial in $n$ variables, and let $\hcal$ be a generic smooth fiber of $W$. In \cite{RSTZ} Ruddat-Sibilla-Treumann-Zaslow give a combinatorial recipe for a skeleton for $\hcal$. In this paper, we show that for a suitable exact symplectic structure  on $\hcal$, the RSTZ-skeleton  can be realized as the Liouville Lagrangian skeleton. 
\end{abstract}

\maketitle
 
Let $(M, \omega = d\lambda)$ be an exact symplectic manifold, and let  $X=X_\lambda$ be  the Liouville vector field defined by $\iota_X \omega = -\lambda$. If $(M, \omega, \lambda, X)$ is a {\em Liouville manifold} (see \cite[Chapter 11]{CE} for definition), then $X$ shrinks $M$ to a compact isotropic (possibly singular) submanifold $\Lambda$, called the {\em Liouville skeleton}. 
The Liouville skeleton is useful for sympletic topology, since  the tubular neighborhood of the skeleton is symplectomorphic to the original manifold up to rescaling the symplectic form. 

A large class of Liouville manifolds come from Stein manifolds, e.g. affine hypersurfaces $\hcal$ in $(\C^*)^n$.  Given an exhausting psh function $\varphi$ on the Stein manifold,  we can define the Liouville structure by setting $\omega = -dd^c \varphi$ and $\lambda = - d^c \varphi$. In \cite{RSTZ}, Ruddat-Sibilla-Treumann-Zaslow give a combinatorial recipe for a topological skeleton in affine hypersurfaces. The RSTZ-skeleton depends on the Newton polytope $Q$ of the defining polynomial for the hypersurface and a star triangulation of $Q$.

It is conjectured that the RSTZ-skeleton can be realized as a Liouville skeleton for a suitable choice of Liouville structure on the hypersurface. Here we construct such Liouville structure using tropicalization. The main idea is contained in the following example: 
\subsection{Example: the pair-of-pants.}
Consider the hypersurface  
\[ \hcal=\{x+y=1\}, \quad x,y \in \C^*. \] The hypersurface can be identified as $\C \RM \{0, 1\}$, a 'pair-of-pants'. A topological skeleton can be constructed as following: fix an arbitrarily small positive number $\epsilon$, and define the skeleton as
\[ \Lambda = \left( \{ |x|=\epsilon \} \cup \{|y|=\epsilon \} \cup \{x\geq \epsilon, y\geq \epsilon\}\right) \bigcap \{x+y=1\}. \]
Thus $\Lambda$ has the shape of two circles connected by an interval. To realize it as a Lagrangian skeleton, we need to choose an exact symplectic structure. Consider the following function $\varphi$ on $(\C^*)^2$ and its restriction on the pair-of-pants
\[ \varphi(x,y) = (\log |x|  - \log \epsilon)^2 + (\log |y| - \log \epsilon)^2. \]
It is easy to check  that $\varphi$ is a psh function on $(\C^*)^2$, and restricts to be a psh function on any complex submanifold of $(\C^*)^2$. Geometrically, $\varphi$ is constructed by taking the projection map
\[ \Log=\log | \cdot |: (\C^*)^2 \to \R^2 \]
and then taking Euclidean distance on $\R^2$ to a point 
\[ \varphi(z) = |\Log(z) - p_\epsilon|^2, \quad  p_\epsilon=(\log \epsilon, \log \epsilon). \] 
The hypersurface $\{ x+y=1\}$ projects under $\log | \cdot |$ to an 'amoeba' shaped region in $\R^2$, with three tenacles asymptotic to the three rays (drawn as dashed lines in the figure)
\[ \{ \log|x|=0, \log|y| \ll 0\}, \quad \{ \log|x|\ll 0, \log|y| = 0\}, \quad \{ \log|x|= \log|y| \gg 0\}. \]
The Louville flow $X_\lambda$ on $\hcal$ induced by $\lambda = - d^c\varphi$ is the same as the negative gradient flow $-\nabla_\omega(\varphi)$ of $\varphi$ with respect to the Kahler metric  $\omega = -dd^c \varphi$. The critical point of $\varphi$ on $\hcal$ can be identified with the critical points on the amoeba $\Log(\hcal)$ with respect to the distance function to point $p_\epsilon$. The unstable manifold of $-\nabla \varphi$ on $\hcal$ is topologically two circles together with an interval. 
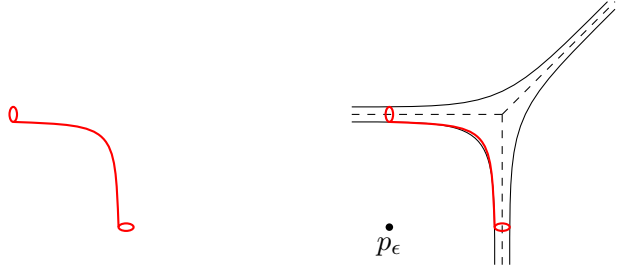
\begin{figure}
\begin{tikzpicture}
\begin{scope}[xshift=-5cm]
\draw [thick, red] (-1.5, 0) ellipse (0.05 and 0.1);
\draw [thick, red] (0, -1.5) ellipse (0.1 and 0.05);
\draw [thick, red](-1.5, -0.1) .. controls (-0.15, -0.15) .. (-0.1, -1.5); 
\end{scope}
\draw [dashed] (0,0) -- (1.5,1.5); 
\draw [dashed] (0,0) -- (-2,0); 
\draw [dashed] (0,0) -- (0,-2); 
\draw (-2,0.1) .. controls (0, 0.1) .. (1.4, 1.5); 
\draw (0.1, -2) .. controls (0.1, 0 ) .. (1.5, 1.4); 
\draw (-2, -0.1) .. controls (-0.1, -0.1) .. (-0.1, -2); 
\fill (-1.5, -1.5) circle (0.05) node [anchor=north] {$p_\epsilon$}; 
\draw [thick, red] (-1.5, 0) ellipse (0.05 and 0.1);
\draw [thick, red] (0, -1.5) ellipse (0.1 and 0.05);
\draw [thick, red](-1.5, -0.1) .. controls (-0.15, -0.15) .. (-0.1, -1.5); 
\end{tikzpicture}
\caption{RSTZ-skeleton (on the left) and its embedding for the pair-of-pants. \label{pop}}
\end{figure}

\subsection{Set-up and Summary of results}
To state our main result precisely, we need some notation. 

Let $M, N$ be dual lattices of rank $n$. Let $T = \R / 2 \pi \Z$. For any abelian group $G$, e.g. $G = \C^*, \R, T$, we define $M_G := M \ot_\Z G$ and similarly for $N_G$. If we fix a basis of $M$, then $M \cong \Z^n$, and $M_{\C^*}\cong (\C^*)^n, M_\R \cong \R^n, M_T \cong T^n$. 

Let $Q \subset N_\R$ be a integral convex polytope of full-dimension containing $0$.\footnote{The case where $Q$ is not full-dimension can be reduced to this one, by defining $N'_\R =\spann_\R(Q) \subset N_\R$, and $M_\R \onto M'_\R$. The skeleton for $(M_\R, N_\R, Q)$ would be that of $(M'_\R, N'_\R, Q)$ times $(S^k)$ where $d=\dim M_\R - \dim M'_\R$. }
  Let $\tcal$ be a coherent star triangulation of $Q$ based at $0$, and let $\pa \tcal$ be the subset of $\tcal$ with simplices not containing $0$ as a vertex. Let $\Sigma_\tcal$ be the simplicical fan spanned by the simplices in $\tcal$. Let $A$ denote the vertices of $\tcal$, and $\pa A$ that of $\pa \tcal$, so that $A = \pa A \cup \{0\}$. 

We fix two functions
\[ h: A \to \R, \quad \Theta: A \to T, \]
 such that $h$ induces the coherent star triangulation of $\tcal$. Without loss of generality, we let 
 \[ h(0) = 0, \quad \Theta(0)=\pi. \]
Let $\wh h: Q \to \R$ denote the convex piecewise linear function on $Q$ extending $h$. 


We define conical Lagrangian $\Lambda_{\tcal, \Theta} \subset M_T \times N_\R \cong T^* M_T$ by
\be \Lambda_{\tcal, \Theta} := \bigcup_{\tau \in \pa \tcal} \{ \theta \in M_T: \la \alpha,  \theta \ra = \Theta(\alpha) \text{ for all vertices $\alpha \in \tau$} \} \times \cone(\tau) \label{e:intro-lambda}\ee
where we used the pairing $\la -, - \ra: M_T \times N_\R \to T$ induced by the canonical pairing between $M, N$, and $\cone(\tau) = \R_{\geq 0} \times \tau$. 
We also define the {\em generalized RSTZ-skeleton} \cite{RSTZ} by
\be \label{RSTZ-skeleton}
 \Lambda_{\tcal, \Theta}^\infty := \bigcup_{\tau \in \pa \tcal} \{ \theta \in M_T: \la \alpha,  \theta \ra = \Theta(\alpha) \text{ for all vertices $\alpha \in \tau$} \} \times \tau \ee
 
\brem
The definition of the original RSTZ-skeleton is for $\Theta|_{\pa \acal}=0$ and is living over the boundary of the Newton polytop $|\pa \tcal| \subset \pa Q \subset N_\R$. In this case $\Lambda_{\tcal, \Theta}$ is first defined by \cite{FLTZ1}.  We will sometimes identify $|\pa \tcal|$ with its projection to $N_\R^\infty :=(N_\R \RM \{0\})/ \R_{>0}$, then and $\Lambda_{\tcal, \Theta}^\infty$ is homeomorphic to the 'end-at-infinity' of  $\Lambda_{\tcal, \Theta}$, hence the notation. 
\erem

For all large enough $\beta > 0$, we define the {\em tropical  polynomial} as
\be f_{\beta, h, \Theta}(z) = \sum_{\alpha \in  A} e^{-i \Theta(\alpha)} e^{-\beta h(\alpha)} z^\alpha. \label{d:trop-f}\ee
where $z^\alpha$ is a monomial   on $M_{\C^*} \cong (\C^*)^n$. Let
\[ \hcal_{\beta, h, \Theta}:= \{ z \in M_{\C^*} \mid  f_{R, h, \Theta}(z) =0 \} \]
denote the complex hypersurface defined by $f_{R, h, \Theta}$.

\btu [\cite{RSTZ}]
If $\Theta|_{\pa \acal}=0$, then the skeleton $\Lambda^\infty_{\tcal, \Theta}$ embeds into the hypersurface $\hcal_{R, h, \Theta}$ as a strong deformation retract. 
\etu

We prove the following theorem, for general $\Theta$. 
\begin{maintheo}
The hypersurface $\hcal_{R, h, \Theta}$ admits a Liouville structure such that its Liouville skeleton is homeomorphic to $\Lambda^\infty_{\tcal, \Theta}$. 
\end{maintheo}

\brem
The deformation of $ f_{\beta, h, \Theta}$ by varying $\Theta$ continuously induces equivalences of categories for Fukaya-Seidel category $FS(M_\CS, f_{\beta, h, \Theta})$. \footnote{We thank Gabe Kerr for the suggestion to consider this freedom of coefficients.}  By realizing the dependence on $\Theta$ in Lagrangian skeleto $\Lambda_{\tcal, \Theta}$, we can prove \cite{Z} equivalences of categories among the infinitesimally wrapped Fukaya categories $\z{Fuk}(T^*M_T, \Lambda_{\tcal, \Theta})$ used in \cite{FLTZ1, FLTZ2}. 
\erem

\subsection{Sketch of Proof}
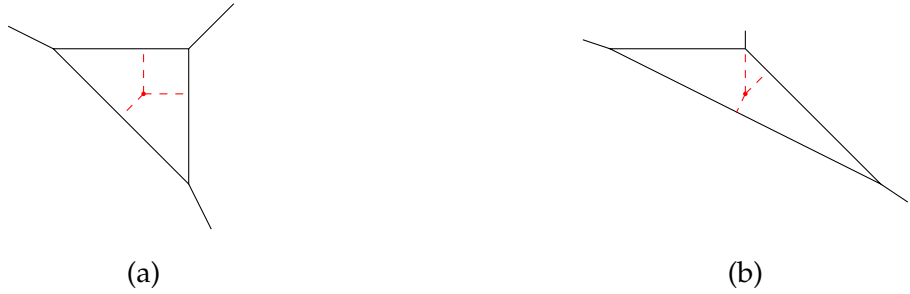
\begin{figure}[h] 
\begin{tikzpicture}
\begin{scope} [xshift=-4cm, scale=0.6]
\draw (1,1) -- (1,-2) -- (-2, 1) -- cycle; 
\draw (1,1) -- (2,2); 
\draw (1,-2) -- (1.5, -3); 
\draw (-2, 1) -- (-3, 1.5); 
\fill [red] (0,0) circle (0.05); 
\draw [red, dashed] (0,0) -- (1,0); 
\draw [red, dashed] (0,0) -- (0,1); 
\draw [red, dashed] (0,0) -- (-0.5,-0.5); 
\draw (0, -4) node {(a)};  
\end{scope}

\begin{scope} [xshift=+4cm, scale=0.6]
\draw (0,1) -- (3,-2) -- (-3, 1) -- cycle; 
\draw (0,1) -- (0,1.4); 
\draw (3,-2) -- (3.6, -2.4); 
\draw (-3, 1) -- (-3.6,1.2); 
\fill [red] (0,0) circle (0.05); 
\draw [red, dashed] (0,0) -- (0.5,0.5); 
\draw [red, dashed] (0,0) -- (0,1); 
\draw [red, dashed] (0,0) -- (-0.2,-0.42); 
\draw (0, -4) node {(b)}; 
\end{scope}
\end{tikzpicture}
\caption{\label{f:adapt} (a) Tropical amoeba of $f=1+e^{-\beta}(x+y+1/xy)$ (b) Tropical amoeba of $f=1+e^{-\beta}(x+xy+1/xy^2)$. The function $|x|^2$ on $\R^2$ is  adapted to the amoeba polytope in (a) but not the one in (b), since the minimum of $|x|^2$ on the top edge lies on the endpoint.  }
\end{figure}

The idea of the proof is illstrated in the above example, that is,  we project the hypersurface $\hcal_{\beta, h, \Theta}$ to $M_\R$, then use a distance function to a point to induce the psh function $\varphi$, which in turn induces a Liouville structure on $\hcal$. However, there are  two technical modification used. 
\begin{enumerate}
\item The first modification is is to 'straighten the tube', or called tropical localization by Mikhalkin and Abouzaid \cite{Mi, Ab}. In the defining Laurent polynomial $f=f_{\beta,h,\Theta}$, not all terms are of equal importance at all points on $\hcal = \hcal_{\beta,h,\Theta}$. We may drop the irrelevant terms and simplify the  defining equation for $\hcal$ locally. \bigskip
\item The second modification is to find a convex function $\varphi$ on $M_\R \cong \R^n$ adapted to the 'tropical amoeba' of the hypersurface. Let $\ccal_0=:P$ is denote the convex polytope for the complement of the tropical amoeba corresponding to vertex $0$. The condition for $\varphi$ is (Definition \ref{d:adapt}) (1) $\varphi(\lambda x) = \lambda^2 \varphi(x)$ for $\lambda > 0, x \in M_\R \RM\{0\}$, and (2) for each face $F$ of $P$ of positive dimension, we want $\varphi|_F$ to have a minimum in $\Int(F)$. See Figure \ref{f:adapt} for an example, where $|x|^2$ fails to be a good potential for case (b).  
\end{enumerate}

\brem
The use of 'non-standard' Kahler potential $\varphi$ on $M_\CS$ (non-canonically isomorphic to  $(\C^*)^n$) may be unorthodoxical, but it is natural in some sense. (1) The often used 'standard' Kahler potential $\sum_i (\log |z_i|)^2$ on $M_\CS$ is not standard in the first place, since it depends of the choice of basis for $M$.  (2) To identify the RSTZ-skeleton \eqref{RSTZ-skeleton} that lives in $M_T \times N_\R$ with the Liouville skeleton \eqref{Liouville-skel} that lives in $M_\CS \cong M_T \times M_\R$, one needs to identify $N_\R$ with $M_\R$.  Here this is done using the Legendre transformation induced by $\varphi$. Equivalently, one fixes a (Finsler) metric on fibers of $TM_T$ and identify $TM_T \cong T^*M_T$. 
\erem

\subsection{Related works}
The study of skeleton for Liouville (or Weinstein) manifold was motivated largely by Homological Mirror Symmetry (HMS).  It was Kontsevich's original proposal, to compute the Fukaya category of a Weinstein manifold $W$ by taking global section of a (co)sheaf of categories living on the skeleton. Following this approach, the (complex) 1-dimensional case has been studied by  \cite{STZ, PS} and  \cite{DK}. 
In general, the category under consideration has two versions, a microlocal sheaf version, and a Floer-theoretic version. The two versions are expected to agree by the on-going work of Ganatra-Pardon-Shende \cite{GPS1, GPS2}.  The microlocal sheaf version, originates from the seminal work of Nadler and Zaslow \cite{NZ, N-brane}, says the infinitesimally wrapped Fukaya category on $T^*M$ with asymptotic condition of non-compact Lagrangian given by a conical Lagrangian $\Lambda$, is equivalent to constructible sheaves on $M$ with singular support in $\Lambda$. The wrapped Fukaya category also has a microlocal sheaf version, developed by Nadler \cite{N4}. 

The microlocal sheaf category for local Lagrangian singularities has been studied by Nadler. In \cite{N-arb}, Nadler defined a class of 'simple' singularities, termed 'arboreal singularities', and proved that the microlocal sheaf category on arboreal singularity is equivalent to the category of representation of quivers. In \cite{N-leg}, Nadler showed that one can deform an arbitrary Lagrangian singularity to an arboreal one, while keeping the microlocal sheaf category invariant. It is also expected that such arborealization can be induced by a perturbation of Weinstein structure\cite{St, ENS}. 

The skeleton for $n$ dimensional pair-of-pants $\pcal_n$ has been studied by Nadler\cite{N4}, where a higher dimensional analog of Figure \ref{pop} is constructued. A $\Sigma_{n+2}$-symmetric skeleton for $\pcal_n$ is constructed by Gammage-Nadler \cite{GN}, where $\Sigma_{n+2}$ is the symmetric group. With the technique of tropical phase variety of Kerr-Zharkov \cite{KeZh}, we hope to find other Lagrangian skeleta $\Lambda_{k}$ for $\pcal_{n-1}$, where $k=1, \cdots,n$ indicating the number of dominant terms in the defining equation of $\pcal_{n-1}$.  

In Gammage-Shende \cite{GS}, as one ingredient in proving HMS for the toric boundary of a toric variety, they constructed the Liouville skeleton for the same hypersurface as considered here. However, their results \cite[Theorem 3.4.2]{GS} depends on the following hypothesis that, {\em there exists some tropicalization function $h: A \to \R$ and some identification $M \cong \Z^n$, such that the tropical amoeba polytope $P=\{x \in M_\R: \la x, \alpha \ra \leq h(\alpha), \forall \alpha \in \pa A\}$ contains $0$ as an interior point, and $|x|^2$ restricts to each face $F$ of $P$ has a minimum in the interior of $F$.} This hypothesis is true in two-dimension, and can be verified in certain examples, e.g. mirror to weighted projective spaces. But in general, one does not know if it is always true, it would be interesting to find a proof or construct a counter-example. \footnote{We thank Gammage and Shende for this clarification.}
Our approach here does not rely on this hypothesis, which is equivalent to {\em "$|x|^2$ is adapted to the tropical amoeba polytope $P$" for some choice of $h$}. We avoid this by considering a more flexible choice of Kahler potentials on  than $|x|^2$, and our approach works for any choice of $h$ compatible with $\tcal$.

\subsection{Outline} In section \ref{s:trop}, we reviewed the tropical localization of Mikhalkin and Abouzaid. We are careful in picking the cut-off functions such that the inner boundary of the amoeba remains convex (Section \ref{ss:convex}). In Section \ref{s:kahler}, we review how to identify $M_\CS$ with $T^* M_T$ by choosing a Kahler potential, and we introduce the key concept of {\bf Kahler potential adapted to a polytope} in Section \ref{s:smart-pot}. Then we state our main theorems in more detail in Section \ref{main}.

\subsection{Acknowledgements}
I would like to thank my advisor Eric Zaslow for introducing the idea of Lagrangian skeleton and the problem of finding Lagrangian embedding. The idea of using Legendre transformation to identify $N_\R$ and $M_\R$ was inspired by a talk of Helge Ruddat.  I thank Vivek Shende and Ben Gammage for the clarification of their work. I also thank Nicol\`o  Sibilla,  Gabe Kerr, David Nadler and Ilia Zharkov for many helpful discussions.

\section{Tropical Geometry\label{s:trop} } 

\subsection{Triangulation and Amoeba}
We follow \cite[Chapter 7]{GKZ} and \cite{Mi} to give background on coherent triangulations and tropical amoeba. 

Let $A \subset N \cong \Z^n$ and $Q=\conv(A)$ its convex hull. Assume $Q$ has full dimension.  A {\em coherent triangulation} $(\tcal, \psi)$ for pair $(Q, A)$ is a triangulation $\tcal$ of $Q$ with vertices  in $A$ and a piecewise linear (PL) convex function $\psi: Q \to \R$, such that the maximal linear domains of $\psi$ are exactly the maximal simplices of $\tcal$.  
For any assignment $h: A \to \R$, there exists a maximal PL convex function $\h h: Q \to \R$, such that $\h h(\alpha) \leq h(\alpha)$ for all $\alpha \in A$. For generic choice of $h$, $\h h$ induces a triangulation $\tcal$ of $(Q, A)$. 
Since we work with a fixed triangulation instead of considering all possible triangulations, we will assume $A$ is the set of vertices of $\tcal$ by reducing $A$ if necessary. If $(\tcal,\psi)$ is a coherent triangulation of $(Q,A)$, $h=\psi|_A$, we may also denote $(\tcal, \psi)$ by $(\tcal, h)$.

Let $(\tcal, \psi)$ be a coherent triangulation of $(Q,A)$. We define the Legendre transformation of $\psi$ as
\[ L_\psi: M_\R \to \R, \quad  u_\psi(y) = \max_{x \in A} \la x, y \ra - \psi(y),  \]
where $\la -, - \ra$ is the dual pairing $M_\R \times N_\R \to \R$. 
One can show that $L_\psi$ is a PL convex function on $M_\R$, inducing a cell-decomposition of $M_\R$ dual to the triangulation of $\tcal$ on $Q$. If $\tau \in \tcal$ is a $k$-simplex, then we use $C_\tau$ or $\tau^\vee$  to denote the dual cell of dimension $n-k$ in $M_\R$. In particular, $C_\alpha$ are the $n$-dimensional cells of $M_\R$.  The cells and simplices are closed in our convention. 

\bd
The {\em tropical amoeba} $\Pi_\psi \subset M_\R$ is defined as the singular loci of $L_\psi$. 
\ed

The tropical amoeba is the limit of amoeba, which we now define. Given a coherent triangulation $(T,h)$ of $(Q,A)$, $h: A \to \R$,  we may define the patchworking polynomial
\[ f_{\beta, h}(z) = \sum_{\alpha \in A} e^{-\beta h(\alpha)} z^\alpha: M_\C^* \to \C. \]
More generally, given a function $\Theta: A \to T$, we have 
\[ f_{\beta, h, \Theta} (z) = \sum_{\alpha \in A} e^{-\beta h(\alpha)} e^{-i \Theta(\alpha)} z^\alpha: M_\C^* \to \C. \]

\bd
The {\em Log amoeba} $\Pi_{\beta, h, \Theta}$ of $f=f_{\beta, h, \Theta}$ is defined as the image of $f^{-1}(0)$ under the (rescaled) logarithm map
\[ \Log_\beta: M \ot_\Z \C^* \to M \ot_\Z \R, \quad m \ot z \mapsto m \ot \beta^{-1} \log|z|. \]
\ed

Mikhalkin proved the following convergence theorem: 
\btu[\cite{Mi}]
The tropical amoeba $\Pi_h$ is the Hausdorff-limit of rescaled amoeba $\Pi_{\beta,h,\Theta}$ as $\beta \to \infty$. 
\etu

\subsection{Monomial cut-off functions} 
The complements of the tropical amoeba has a one-to-one correspondence with the vertices of the triangulation $\tcal$,
\[ M_\R \RM \Pi_h = \bigsqcup_{\alpha \in A} C_\alpha, \quad M_\R \RM \Pi_{\beta,h} = \bigsqcup_{\alpha \in A} C_{\alpha, \beta}. \]
$C_\alpha$ are convex polyhedra, and $C_{\alpha, \beta}$ are smooth strictly convex domains \cite[Chapter 6, Cor 1.6]{GKZ}. 

The idea of introducing a monomial cut-off function $\chi_{\alpha,\beta}(z)$ is to turn off the term $e^{-\beta h(\alpha)} z^\alpha$ if it is much smaller than the rest, thus straighting the hypersurface. The idea is first used in Abouzaid \cite{Ab} to control the symplectic geometry of the hypersurface.

\begin{figure}[h]
\includegraphics[width =0.4 \textwidth]{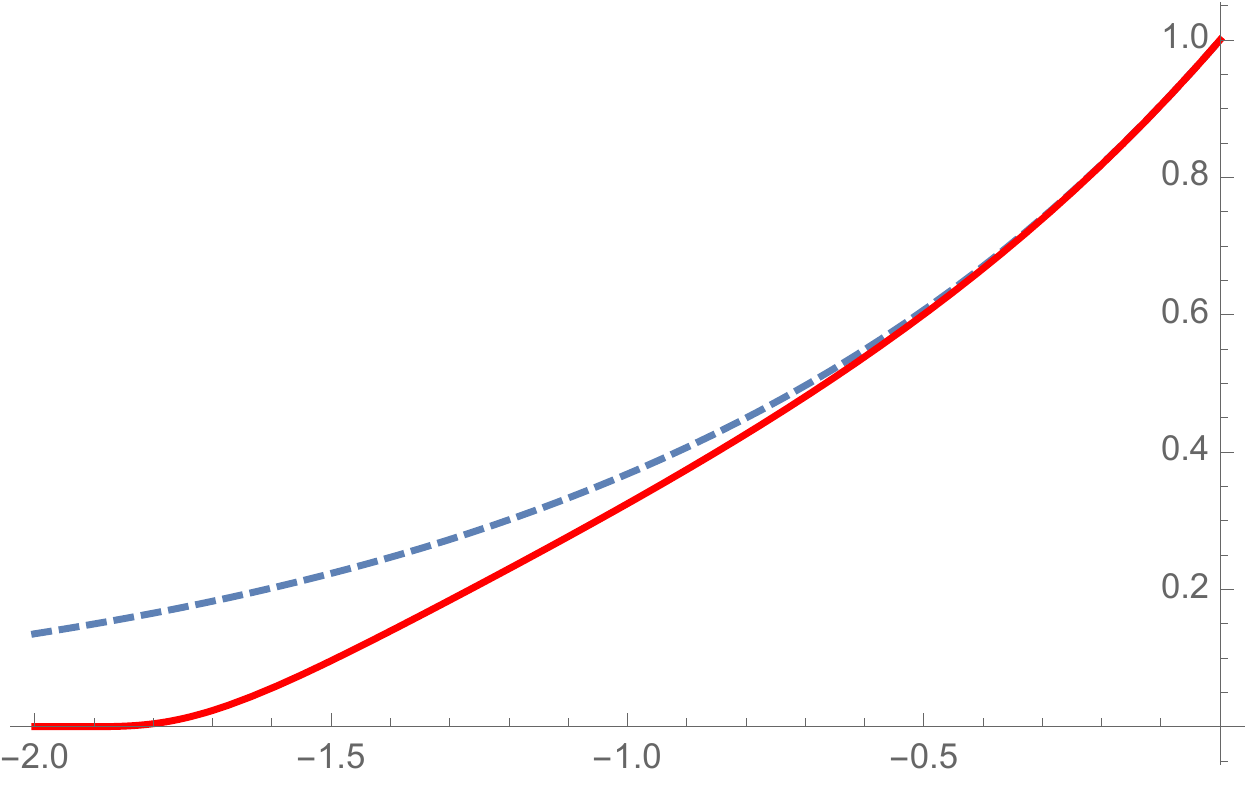}
\caption{\label{chi}. We modify exponential function $e^x$ (dashed line)  to $e^x \chi(x)$(solid line), such that $e^x \chi(x)$ remains convex. }
\end{figure}

We fix a cut-off function $\chi(x)$ on $\R$ with the following properties
\begin{itemize}
\item $\chi(x) = \bcs 1 & x \in [0,\infty) \\ \in (0,1) & x \in (-2,0) \\ 0& x \in (-\infty, -2] \ecs $. 
\bigskip
\item $\chi(x) \exp(x)$ is convex. 
\end{itemize}
\bex
One can check that the following specification of $\chi(x)$ on $[-2,0]$ gives a $C^2$ function on $\R$ with the desired convexity. 
\[ \chi(x) =  e^{-1/(x+2) + 1/2 -x/4+x^2/8}. \]
See Figure \ref{chi} for a plot of $\chi(x) e^{x}$. We have 
\[ \frac{(e^x \chi(x))''}{e^x \chi(x)} = 256 + 608 x + 576 x^2 + 288 x^3 + 85 x^4 + 14 x^5 + x^6, \quad x \in (-2,0). \]
which can be verified to be positive on $(-2,0)$ The non-smooth point of $\chi(x)$ is at $x=0$, and can be mollified if needed. 
\eex

We use log coordinates $(\rho, \theta) \in M_\R \times M_T$ for point $z \in M_\CS$. We also use $\beta$-rescaled log coordinates $(u, \theta) = (\beta^{-1} \rho, \theta)$, thus $u = \Log_\beta(z)$. For each $\alpha \in A$, define  linear function on $M_\R$
\[ l_\alpha(u) := \la u, \alpha \ra - h(\alpha). \]
\bd\label{d:cutoff}
For any vertex $\alpha \in A$, we define the {\em monomial cut-off function} as
\[ \chi_{\alpha, \beta}(u) = \prod_{\text{$\alpha'$ adjacent to $\alpha$ in $\tcal$}} \chi_{\alpha, \alpha', \beta}(u)\] where \[  \chi_{\alpha, \alpha', \beta}(u)=\chi(\beta(l_{\alpha}(u)  - l_{\alpha'}(u)) + \sqrt{\beta}). \]
\ed

We define a distance-like function to region $C_\alpha$,
\[ r_\alpha(u) := L_{\h h}(u) - l_{\alpha}(u)) = \max_{\alpha' \in A} (l_{\alpha'} (u) - l_\alpha(u)).\] 
Thus $r_\alpha(u)$ is a non-negative PL convex function, vanishes only on $C_\alpha$.  

\bp
For all large enough $\beta$, $\chi_{\alpha, \beta}(u)$ satisfies the following property
\[ \chi_{\alpha, \beta}(u) = \bcs 1 & r_\alpha(u) < \beta^{-1/2} \\
0 & r_\alpha(u) >  \beta^{-1/2} + 2 \beta^{-1}
\ecs
\]
\ep
\bpf
If $r_\alpha(u) < \beta^{-1/2}$, then for all $\alpha'$ adjacent to $\alpha$, we have 
\[ l_{\alpha'}(u) - l_{\alpha}(u) < \beta^{-1/2} \quad \RA \quad \beta(l_{\alpha}(u)  - l_{\alpha'}(u)) + \sqrt{\beta}) > 0 \]
thus each factor in $\chi_{\alpha, \beta}(u)$ equals $1$. The other case is similar to check, where $\beta$ large enough means $r_\alpha(u)^{-1}(\beta^{-1/2} + 2 \beta^{-1})$ intersects all the neighboring cells $C_{\alpha'}$ for $C_\alpha$. 
\epf

\bd \label{d:bad-region}
For each $\alpha \in A$, we define the {\em bad region} as the open set 
\[ B_{\alpha,\beta} = \{u \in M_\R \mid \beta^{-1/2} + 2 \beta^{-1} >  r_\alpha(u) >\beta^{-1/2}  \}. \]
The {\em (total) bad region} $B_\beta$ is defined as the union of all $B_{\alpha,\beta}$. The {\em good region} is defined as the closed set $G_\beta := M_\R \RM B_\beta$.  
\ed

On the good region, each $\chi_{\alpha, \beta}$ is either $0$ or $1$, hence we have a partition labeled by cells of $\tcal$: 
\[ G_{\beta} = \bigsqcup_{\tau \in \tcal} G_{\beta, \tau} ,\]
where $G_{\beta, \tau} = \{u \in G_{\beta} \mid \chi_{\alpha, \beta}(u) = 1 \iff \alpha \in \tau \}$ is a closed convex polyhedron with non-empty interior.  

\ss{Tropical Localized Hypersurfaces}
Following \cite[Section 4]{Ab}, we define a family of hypersurfaces $\hcal_s$ as the zero-loci of $f_s(z)$: 
\[ f_{s}(z):= \sum_{\alpha \in A} f_{\alpha} (z) \chi_{\alpha,s}(z), \quad \hcal_s = f_s^{-1}(0), \]
where 
\[  f_\alpha = e^{-\beta h(\alpha) - i \Theta(\alpha)} z^\alpha, \quad \chi_{\alpha,s}(z) = s \chi_{\alpha}(z) + (1-s).   \]
Here we have dropped the $\beta, h, \cdots$ subscripts from previous notations for clarity. These hypersurfaces $\hcal_s$
 interpolates between the complex hypersurface $\hcal$ and tropical localized hypersurface $\wt \hcal$, where
 \[ \hcal:=\hcal_0 , \quad \wt \hcal := \hcal_1. \]


The following proposition is a modification of Proposition 4.2 in \cite{Ab}. 
\bp \label{p:symptrop}
Fix any identification $M_\CS \cong \CSN$. Let $\omega$ be any toric Kahler metric on $M_\CS$, comparable with $\omega_{0} =  i \sum_i d \log z_i \wedge d \wb{\log z_i}$. Then, for all large enough $\beta$, the family of hypersurfaces $\hcal_s$  are symplectic with induced symplectic form from $\omega$.
\ep
\bpf
We proceed as in Proposition 4.2 in \cite{Ab}: to prove $f^{-1}(0)$ is symplectic, it suffices to prove $|\dbar f(z)|_\omega < |\d f(z)|_\omega$. Since $\omega$ and $\omega_0$ are comparable, in the following the norm on $M_\R$ or its dual $N_\R$ are taken to be the Euclidean norm. 

Assume $s>0$ and $z$ is in the 'bad region'(see Definition \ref{d:bad-region}) of the hypersurface $\hcal_{\alpha,s}$, since otherwise the hypersurface is holomorphic at $z$ and there is nothing to prove. Let $I(z)=\{\alpha_0, \cdots, \alpha_k\} \subset A$ be subset of vertices where $r_{\alpha_i}(z) \leq \beta^{-1/2}$, thus $I'(z)$ is vertex set for a simplex $\tau(z)$ of $\tcal$. Since $f_{\alpha_0}(z)^{-1} f_s(z)$ is an equally good defining equation for $\hcal_s$, hence without loss of generality, we may assume that $\alpha_0=0$. Let $I(z) = I'(z) \RM \{\alpha_0\}$. 

Define
\[ F(z) = \sum_{\alpha \in A} |f_\alpha(z)| = \sum_{\alpha in A} e^{\beta l_\alpha(u)} =: e^{\beta \varphi_{\beta,h}(u)}. \]
Then we have $\varphi_{\beta, h}(u) \geq L_h(u)$ for all $u \in M_\R$, and as $\beta \to \infty$ we have $\varphi_{\beta, h}(u) \to L_h(u)$ in $C^0$. 

We first note that of derivatives for the cut-off functions $\chi_{\alpha}(z)$ have a uniform bound
\bea |d \chi_\alpha(\rho)| &=& | d( \prod_{\alpha'} \chi(\sqrt{\beta}+l_\alpha(\rho) - l_{\alpha'}(\rho))) | \\
&\leq& \sum_{\alpha'} | d(\chi(\sqrt{\beta}+l_\alpha(\rho) - l_{\alpha'}(\rho)))| \\
&\leq& \|\chi'\|_{\infty} \sum_{\alpha'} |\alpha'-\alpha| < C
\eea
where the product or sum are over $\alpha'$ adjacent to $\alpha$ in triangulation $\tcal$, and we used the bound $\chi \leq 1$.

We have for $\dbar f_s(z)$. 
\bea
|F(z)^{-1}\dbar f_s(z)| &=& |F(z)^{-1}\sum_\alpha f_\alpha(z) \dbar \chi_{\alpha,s}(z)| 
\leq   s\sum_\alpha e^{\beta (l_\alpha(u) -  \varphi_{\beta,h}(u))} |\dbar \chi_\alpha(z)| \\
&\leq & \sum_\alpha e^{\beta (l_\alpha(u) -  L_h(u))} |\dbar \chi_\alpha(z)| 
 \leq \sum_\alpha e^{\beta (-\beta^{-1/2})} |\dbar \chi_\alpha(z)| \\
&\leq& C e^{-\sqrt{\beta}}
\eea
where we used $\varphi_{\beta,h}(u) \geq L_{h}(u)$ and on the support of $d \chi_\alpha(z)$ we have $r_\alpha(z) = L_h(u) - l_\alpha(u)>1/\sqrt{\beta}$. 

Next, we compute $\d f_s(z)$. 
\bea
|F(z)^{-1}\d f_s(z)| &=& |F(z)^{-1}\sum_\alpha \d f_\alpha(z) \chi_{\alpha,s}(z)+ s f_\alpha(z) \d \chi_\alpha(z)| 
\\&=& F(z)^{-1} |\sum_{i=1}^k \d f_{\alpha_i}(z)| +  O(e^{-\sqrt{\beta}}) \\
&=& F(z)^{-1} |\sum_{i=1}^k f_{\alpha_i}(z) \la \alpha_i, d (\rho + i \theta)| +  O(e^{-\sqrt{\beta}}) \\
&=&   F(z)^{-1}  \left[ \sum_{i=1}^k \sum_{j=1}^k   f_{\alpha_i} (z) \wb{f_{\alpha_j}(z)} \la \alpha_i, \alpha_j \ra \right]^{1/2} +  O(e^{-\sqrt{\beta}}) \\
&>& C_1 F(z)^{-1} \left[\sum_{i=1}^k |f_{\alpha_i}(z)|^2\right]^{1/2} +  O(e^{-\sqrt{\beta}}) \\
&>& C_2 F(z)^{-1} \sum_{i=1}^k |f_{\alpha_i}(z)|+  O(e^{-\sqrt{\beta}})  = C_2 + O(e^{-\sqrt{\beta}}) 
\eea
where $ O(e^{-\sqrt{\beta}}) $ represent a remainder bounded by $e^{-\sqrt{\beta}}$,  $C_1$ is the smallest eigenvalue of the $k \times k$ matrix $M_{ij}=\la \alpha_i, \alpha_j \ra$, which is non-degenerate since $\{\alpha_0=0, \alpha_1, \cdots\}$ are vertices of a $k$-simplex in $\tcal$. We also used that all the $l^p$ ($1\leq p \leq \infty$) norm on $\R^k$ are equivalent. 

Thus, we have shown $|\d(f_{\alpha_0}(z)^{-1}f_s(z))| > |\dbar(f_{\alpha_0}(z)^{-1}f_s(z))|$ for $z \in f_s^{-1}(0)$, hence $\hcal_s$ are symplectic. 
\epf

\bp[Proposition 4.9 \cite{Ab}] \label{family}
The family of hypersurfaces $\hcal_s$ are symplectomorphic for all $s \in [0,1]$. 
\ep

\ss{Tropical Localization with Convexity \label{ss:convex}}
For log amoeba $\Pi_{h, \beta}$ and tropical amoeba $\Pi_{h}$, their connected components of complement $C_{\alpha, \beta}$ and $C_{\alpha}$ are convex. Let $\wt \Pi_{\beta,h} = \Log_\beta(\wt \hcal)$ be the amoeba of the tropical localized hypersurface, and let $\wt C_{\alpha}$ denote the complements dual to vertex $\alpha \in A$. We would like to show that $\wt C_\alpha$ are close to convex as well. 

\bp
The defining equation for $\wt C_\alpha$ is
\[1 = \sum_{\alpha' \text{ adjacent to } \alpha} e^{\beta(l_{\alpha'}(u) - l_{\alpha}(u))} \chi_{\alpha', \beta}(u)=F_{\alpha}. \]
and $F_{\alpha}$ is convex in the good region, i.e., where all $\chi_{\alpha, \beta}(u)$ are constant with value $0$ or $1$. 
\ep
\bpf
The boundary of the complement $\wt C_\alpha$ is where the dominant term equals the sum of the other non-dominant term. By the tropical localization, there are at most $n$ non-dominant terms for a point $z$ on the boundary (thanks to $\tcal$ being a triangulation). Hence the $\theta_i$ can be chosen, such that the argument of the dominant and non-dominant terms are the same. For the second statement, we note that over the good region, $F_{\alpha}$ is a sum of convex functions. 
\epf

\bd\label{d:convexC}
The convex model $\wh C_{\alpha}$ for $\wt C_{\alpha}$ is defined by $\{u \in  {C_\alpha} \mid \wh F_\alpha (u) = 1\}$, where 
\[ \wh F_{\alpha} = \sum_{\alpha' \text{ adjacent to } \alpha}  e^{\beta(l_{\alpha'}(u) - l_{\alpha}(u))}  \chi_{\alpha', \alpha,\beta}(u) \]
\ed

The two defining functions, $F_\alpha$ and $\wh F_\alpha$, differ by the cut-off functions: $F_\alpha$ uses $\chi_{\alpha',\beta}$ which cuts along the boundary of $C_{\alpha'}$, whereas $\wh F_\alpha$ uses $\chi_{\alpha', \alpha, \beta}$ which cuts along the hyperplane separating $C_\alpha$ and $C_{\alpha'}$. However, on $C_{\alpha}$ the two functions and the hypersurfaces are very close, as the following two propositions show. 

\bp \label{Ck}
For all $k \geq 1$,  there are constant $c_k, c'_k$, such that 
\[ \| F  \|_{C^k( C_{\alpha})} + \| \wh F \|_{C^k( C_{\alpha})} \leq c'_k \beta^k \]
and
\[ \| F - \wh F \|_{C^k( C_{\alpha})} < c_k \beta^k e^{- \sqrt{\beta}}. \]
\ep
\bpf
First, we note that since all $\chi_{\alpha_1, \alpha_2, \beta} \leq 1$, we have $\chi_{\alpha'}(u) < \chi_{\alpha', \alpha}(u)$, thus 
\[  \wh F_\alpha (u) >  F_\alpha (u), \quad \wh C_\alpha \subset \wt C_{\alpha}.    \]
For $u \in \wb C_{\alpha}$ and $\alpha' \text{ adjacent to } \alpha$, if $\chi_{\alpha'}(u) - \chi_{\alpha',\alpha}(u) \neq 0$, then $l_{\alpha'}(u) - l_{\alpha}(u)+\sqrt{\beta} \in (-2,0)$. Hence
\[
\wh F_\alpha (u) - F_\alpha (u) =  \sum_{\alpha' \text{ adjacent to } \alpha} e^{\beta(l_{\alpha'}(u) - l_{\alpha}(u))} (\chi_{\alpha}(u) - \chi_{\alpha',\alpha}(u)) < C e^{-\sqrt{\beta}}. \]
Similarly, taking $k$-th derivative, we have 
\[
|\pa_u^k \wh F_\alpha (u) - \pa_u^k  F_\alpha (u)| < C_k \beta^k e^{-\sqrt{\beta}}. \]
where the norm are taken with respect to Euclidean norm on $\R^n$, after choosing an identification $M_\R \cong \R^n$. This finish the proof. 
\epf

Fix $M_\R \cong \R^n$ and equip $\R^n$ with Euclidean metric. Let $S^* \R^n$ denote the unit cosphere bundle. 
If $C$ is a domain with smooth boundary, we define 
\[ \Lambda_{C} = \{ (p, \xi) \in S^*\R^n \mid p \in \pa C,  \xi \in (T_p \pa C)^\perp \text{ and points outward }\} \]

\bp 
We have the following convergence results: 
\begin{enumerate}
\item In the good region in $C_{\alpha}$, $\pa \wt C_\alpha = \pa \wh C_{\alpha}$. 
\item  The Hausdorff distance between $\pa \wt C_\alpha$ and $\pa \wh C_{\alpha}$ is $O(\beta^{-1} e^{-\sqrt{\beta}})$.  
\item  The Hausdorff distance between $\Lambda_{\wt C_\alpha}$ and $\Lambda_{\wh C_{\alpha}}$ is $O(\beta e^{-\sqrt{\beta}})$. 
\end{enumerate} 
\ep
\bpf
We will write $F=F_\alpha, \wh F = \wh F_\alpha$ and so on, omitting the $\alpha$ subscript when it is not confusing. 

(1) Since in the good region in $C_\alpha$, all the cut-off functions $\chi_{\alpha'}$ and $\chi_{\alpha', \alpha}$ are equal. 

(2) Since $\wh F  \leq F$, hence the domain $\wh C \subset \wt C \subset C$. If $u \in \pa \wh C \RM \pa \wt C$, we take gradient flow of $F$, starting from $u$ and ending on $u' \in  \pa \wt C$.  Let $\gamma:[0,t] \to C_\alpha$ denote this integral curve. Since around $\pa \wt C$ and $\pa \wh C$, we have uniform lower bound for $|dF|$ and $|d\wh F|$ by some constant $c \beta$, hence
\[ c \beta \dist(u, u') \leq c \beta t < \int_0^t |\nabla F(\gamma(s))| ds \]
and 
\[ \int_0^t |\nabla F(\gamma(s))| ds= F(\gamma(t)) - F(\gamma(0)) = 1 - F(u) = \wh F(u) - F(u) < C e^{-\sqrt{\beta}} \]
hence 
\be \label{distu}
 \dist(u, u') = O( \beta^{-1} e^{-\sqrt{\beta}}) \ee
Similarly, if we start from $u' \in \pa \wt C \RM \pa \wh C$ we may find $u \in \wh C$ using gradient flow of $\wh F$, with the same bound as above. This establishes the bound on Hausdorff distance

(3) Let $u \in \pa \wh C \RM \pa \wt C$, and $u' \in  \pa \wt C$ constructed as in (2). We have
\bea 
|d F(u') - d \wh F(u)| \leq |d F(u') - d F(u)| + |d F(u) - d \wh F(u)| \\
\leq  \|F\|_{C^2} \dist(u,u') +O(\beta  e^{-\sqrt{\beta}})   = O(\beta e^{-\sqrt{\beta}})  \eea
where we used the $C^2$ bound of $F$ in Proposition \ref{Ck} and the distance bound in \eqref{distu}. 
\epf

Let $\Lambda_{C_{\alpha}}$ denote the Legendrian of the unit exterior conormal to $C_\alpha$. 

\bp \label{lm:GH-conv} Let $\alpha \in A$, $\wh C_\alpha$ be the tropical localized amoeba's complement. We have the following convergence of the boundary of $\wh C_\alpha$, and its Legendrian lifts: 
\begin{enumerate}
\item  The Hausdorff distance between $\pa \wh C_\alpha$ and $\pa  C_{\alpha}$ is $O(1/ \sqrt{\beta})$.  
\item  The Hausdorff distance between $\Lambda_{\wh C_\alpha}$ and $\Lambda_{C_{\alpha}}$ is $O(1/ \sqrt{\beta})$. 
\end{enumerate} 

\ep
\bpf
Consider the simplices around vertex $\alpha$ in $\tcal$. Let $\tau$ be such a $k$-dimensional simplex, with vertices $\alpha, \alpha_1, \cdots, \alpha_k$. Denote the dual face by $\tau^\vee$ on the polytope $C_\alpha$. 
We also define an locally closed subset $U_\tau \subset \pa \wh C_{\alpha}$,  such that $z \in U_\tau$ iff the set $I_\alpha(z) = \{\alpha_1, \cdots, \alpha_k\}$. 
\[ I_\alpha(z): = \{\alpha' \text{ adjacent to } \alpha,  \chi_{\alpha', \alpha, \beta}(z) > 0 \}.  \]

Define the orthogonal projection projection map 
\[ \pi_\tau: U_\tau \to \tau^\vee. \]
If $u \in U_\tau$, $u' = \pi_\tau(u)$, then since
\[ -1/\sqrt{\beta}- 2/\beta < l_{\alpha_i}(u)-l_{\alpha}(u)  < 0, \quad \text{ and } \quad  l_{\alpha_i}(u')-l_{\alpha}(u') =0 \]
Hence
\[ dist(u, u') < O( \sum_{i=1}^k |l_{\alpha_i}(u)-l_{\alpha}(u)|) = O(1/\sqrt{\beta}). \]
Also $\tau^\vee$ is in $O(1/\sqrt{\beta})$ neighborhood of $Im(\pi_\tau)$, thus the Hausdorff distance between $\tau^\vee$ and $U_\tau$ is $O(1/\sqrt{\beta})$. Considering all faces $\tau^\vee$ of $C_\alpha$ proves the first statement. 

For the second statement, we further note that, for any $u \in U_\tau$, the exterior unit conormal $\xi$ of $\wh C_\alpha$ at $u$ is contained in $\cone(\alpha_1-\alpha, \cdots, \alpha_k-\alpha) = \R_+  \tau$. Define
\[ \Lambda_{\tau^\vee} := \tau^\vee \times (\R_+ \cdot \tau) \; \cap \; S^*\R^n. \]
Then the projection map $\pi_\tau$, lifts to 
\[ \wt \pi_\tau: \Lambda_{\wt C_\alpha}|_{U_\tau} \to \Lambda_{\tau^\vee}, \quad (u, \xi) \mapsto (\pi_\tau(u), \xi) \]
Since the fiber direction has distance zero, the similar argument as (1) proves the second statement.  
\epf

\section{Legendre transformation and Toric Kahler Potential. \label{s:kahler}}
In this section we use Legendre transform to define a diffeomorphism between  $(\C^*)^n$ and $T^*T^n$, and define a Kahler structure on $\CSN$. 

\subsection{Legendre transformation}
Let $V$ be a real vector space of dimension $n$, and $V^\vee$ be its dual space. There is a natural identification of symplectic space
\[ T^*V \cong V \times V^\vee \cong T^* V^\vee. \]
Let $\pi_{V}$ and $\pi_{V^\vee}$ denote the projection of $V \times V^\vee$ to its first and second factor, respectively. 

Let $\varphi$ be a smooth strictly  convex function on $V$. The Legendre transformation for $\varphi$ is defined as
\[ \Phi_{\varphi}: V \to V^\vee, \quad x \mapsto d\varphi(x). \]
 We will always assume $\varphi$ satisfies some growth condition such that the Legendre transformation is surjective.
The Legendre dual $\psi$ of $\varphi$ is also a convex function defined as
\[ \psi: V^\vee \to \R, \quad y \mapsto \sup_{x \in V} \la x, y\ra - \varphi(x) = \la \Phi^{-1}_\varphi(y), y \ra - \varphi(\Phi^{-1}_\varphi(y)). \]

If we fix a linear coordinate $\rho=(\rho_1, \cdots, \rho_n)$ on $V$ and dual coordinate $p=(p_1, \cdots, p_n)$ on $V^\vee$, then the Legendre transformation can be written as
\[ p_i = \pa_{\rho_i} \varphi(\rho). \]
If $p = d \varphi(p)$, then Legendre dual function
\[ \psi(p) = \sum_i \rho_i p_i - \varphi(\rho). \]
And the two matrices $\Hess \varphi(\rho) = \pa_{ij} \varphi(\rho)$ and $\Hess \psi(p) = \pa_{ij} \psi(p)$ are inverse of each other. 
There is a metric on $V$ induced by $\varphi$: 
\[ g_\varphi = \pa_{ij} \varphi(\rho) d \rho_i \ot d \rho_j. \]

The above contruction can be interpreted symplectically. 
Consider the graph Lagrangian $\Gamma_{d \varphi}$ in $T^* V$
\[ \Gamma_{d\varphi} = \{ (x, y) \in V \times V^\vee \mid y = d \varphi(x) \}. \]
Let $L = \Gamma_{d \varphi}$. Then the Legendre transform is $\Phi_\varphi = \pi_{V^\vee}|_L \circ \pi_V|_L^{-1}$
\[ \bcd
& L \arrow[swap]{ld}{\pi_{V}} \arrow{rd}{\pi_{V^\vee}} &\\
V & & V^\vee
\ecd. \]
$L$  as a section in $T^*V^\vee$ is the graph of $\Gamma_{d\psi}$ for the Legendre dual function $\psi$ of $\varphi$.

The following lemma says, gradient vector field on $V$ and differential one form are related by Legendre transformation. 
\bl \label{p:gradf}
Let $\varphi$ be any smooth convex function on $V$, and let $f: V \to \R$ be any smooth function. For any $\rho \in V$,  and $p = \Phi_{\varphi}(\rho) \in V^\vee$,  then $(\Phi_\varphi)_* (\nabla f|_\rho) \in T_p V^\vee \cong V^\vee$ and $df(\rho) \in T^*_\rho V \cong V^\vee$ are equal. 
\el
\bpf
We work with linear coordinates $(\rho_1, \cdots, \rho_n)$ on $V$ and dual coordinate $(p_1, $ $\cdots,$ $ p_n)$ on $V^\vee$. Let $g_{ij} = (g_\varphi)_{ij} = \pa_{ij} \varphi$ and $g^{ij}$ be the matrix inverse of $g_{ij}$. 
\bea 
(\Phi_{\varphi})_* \nabla f(\rho) &=& \sum_{i,j,k} \pa_{\rho_k} f \cdot g^{jk} \cdot \frac{\pa p_i(\rho)}{\pa \rho_j} \cdot \pa_{p_i} \\
&=&   \sum_{i,j,k}  \pa_{\rho_k} f \cdot g^{jk} \cdot g_{ij} \cdot \pa_{p_i} \\
&=&  \sum_{i,k}  \pa_{\rho_k} f \cdot \delta^i_k \cdot \pa_{p_i} = d f.
\eea

\epf

\subsection{Identification between $M_{\C^*}$ and $T^* M_T$.}
\label{ss:liouville}
There is a canonical complex structure on $M_{\C^*} \cong M_\R \times M_T$, and a canonical symplectic structure on $T^* M_T \cong N_\R \times M_T$. We will use notation $\theta \in M_T, \rho \in M_\R$ and $p \in N_\R$.  If we fix a $\Z$-basis for $M$, then we have $M_{\C^*} \cong (\C^*)^n = \{(e^{\rho_i + i \theta_i})_i\}$ and $T^*M_T \cong T^*T^n = \{(\theta_i, p_i)_i\}$.

Let $\varphi: M_\R \to \R$ be a smooth strictly convex function such that the Legendre transformation $\Phi_\varphi: M_\R \to N_\R$ is surjective. We abuse notation and also denote by $\varphi$ the pullback via $M_{\C^*} \to M_\R$, and call $\varphi$ a \Kahler potential on $M_{\C^*}$. Then we may define Liouville one-form and symplectic two-form on $M_{\C^*}$
\[ \lambda = -d^c \varphi, \quad \omega = -dd^c \varphi. \]
In coordinate form, we have 
\[ \lambda_\varphi = \sum_i \pa_i \varphi(\rho) d\theta_i, \quad \omega_\varphi = \sum_{i,j} \pa_{ij} \varphi(\rho) \,d \rho_i \wedge d \theta_j. \]
The Riemannian metric can also be obtained by $g_\varphi(X, Y) = \omega_\varphi(X, JY)$, where $J \pa_{\rho_i} = \pa_{\theta_i}, J \pa_{\theta_i} = - \pa_{\rho_i}$, or in coordinate form
\[ g  = \sum_{i,j} \pa_{ij} \varphi(\rho) (d\rho_i \ot d\rho_j + d\theta_i \ot d \theta_j). \]

If we equip $T^*M_T$ with the standard exact symplectic structure $(\omega, \lambda)$: 
\[ \lambda_{std} = \sum_i p_i d \theta_i, \quad \omega_{std} = \sum_{i} d p_i \wedge d\theta_i, \]
then by Legendre transformation $\Phi_\varphi \times \id: M_{\C^*} =M_\R \times M_T \to N_\R \times M_T = T^*M_T$, we have
\[ (\Phi_\varphi \times \id)^*(\lambda_{std}) = \lambda_{\varphi}, \quad  (\Phi_\varphi \times \id)^*(\omega_{std}) = \omega_{\varphi}. \]

\subsection{Homogeneous \Kahler potential}
Next we will restrict ourselves to homogenous convex functions as \Kahler potential. 
\bd \label{d:homog}
A convex function $\varphi$ on $M_\R$ is said to be {\em homogeneous of degree $d$} for some $d \geq 1$, if for any $0 \neq x \in M_\R$ and any $\lambda > 0$, we have 
\be \varphi( \lambda x) = \lambda^d \varphi(x),  \ee
and $\Omega = \{x : \varphi(x) \leq 1 \}$ is a bounded strictly convex closed set with smooth boundary. 
\ed 

\brem
Any positive definite quadratic form on $M_\R$ is a homogeneous degree two convex function. More generally, any bounded strictly convex subset $\Omega \subset M_\R$ with smooth boundary and containing $0$ as an interior point determines a homogeneous degree $d$ convex function $\varphi_{\Omega,d}$ such that $\Omega = \{x : \varphi(x) \leq 1 \}$. 
\erem

\bp
For any homogeneous convex function $\varphi$ of degree $d$ with $d \in [0,\infty)$,  we have \\
(1) $\varphi$ is smooth on $M_\R \RM \{0\}$. \\
(2) $\varphi$ is $C^k$ at $0$ where $k$ is the largest integer less than $d$. \\
(3) If $d >1$, then $\varphi$ is strictly convex. 
\ep
\bpf
(1) and (3) are easy to verify. We only prove (2). 
Fix a linear coordinate $x_1, \cdots, x_n$ on $M_\R$. For multi-index $\alpha=(\alpha_1, \cdots, \alpha_n)$, any point $0 \neq x \in M_\R$ and $\lambda > 0$, we have 
$ \pa_x^{\alpha} \varphi( \lambda x) = \lambda^{d - |\alpha|} \pa_x^{\alpha} \varphi(x). $ Hence if in addition $|\alpha| \leq k < d$,  then 
$ \lim_{\lambda \to 0} \pa_x^{\alpha} \varphi( \lambda x) = 0.$
Hence all $k$-th derivative can be continuated to $x=0$. 
\epf

\bl
If $\varphi$ is a homogeneous degree $d$ convex function, then for $\lambda>0$
\[ \Phi_\varphi(\lambda \rho) = \lambda^{d-1} \Phi_\varphi(\rho). \]
\el

\bd\label{d:proj-leg}
Let $M_\R^\infty := (M_\R \RM 0)/\R_{>0}$ and $N_\R^\infty := (N_\R \RM 0)/\R_{>0}$. Then we define the projective Legendre transformation
\[ \Phi_\varphi^\infty: M_\R^\infty \to N_\R^\infty. \]
\ed
It is easy to check that $\Phi_\varphi^\infty$ is an orientation perserving diffeomorphism from $S^{n-1}$ to itself. Geometrically, if we take the level set $S = \varphi^{-1}(1)$, then each element in $M_\R^\infty$ corresponds to a point on $S$, and the outward conormal of $S$ at the point is the element in $N_\R^\infty$ obtained by $\Phi_\varphi^\infty$. 

\bp
\label{p:gradphi}
Let $\varphi$ be any homogeneous convex function on $M_\R$ of degree $k>1$, and equip $M_\R$ with metric $g_\varphi$ induced from Hessian of $\varphi$. Then the integral curves in $M_\R \RM \{0\}$ of the gradient of $\varphi$  are rays. Equivalently, 
\[ \nabla \varphi(\rho) = C(\rho) \sum_i \rho_i \pa_{\rho_i}, \quad C(\rho)>0. \]
\ep
\bpf
For any nonzero $\rho \in M_\R$, we have $\Phi_\varphi(\rho) = d \varphi(\rho) $, also by Proposition \ref{p:gradf} we have $ (\Phi_\varphi)_* (\nabla \rho) = d\varphi(\rho) $, hence the gradient vector field   $\nabla \varphi$ on $M_\R$ when pushed-forward to $N_\R$ is exactly the radial vector field $p\pa_p$ whose integral curves are rays. Since $\varphi$ is homogeneous, hence $\Phi_\varphi$ takes ray to ray, hence the integral curve of $\nabla \varphi$ is the pull-back of integral curve of $p \pa_p$, i.e. rays. 
\epf

\subsection{\Kahler potentials Adapted to a Polytope \label{s:smart-pot}}
This is one of the key construction in this paper. We replace the Kahler potential $\sum_i u_i^2$ on $\CSN$ where $u_i = \log|z_i|$ by any homogeneous degree two Kahler potential $\varphi(u)$. 

Let $P$ be a convex polytope (possibly unbounded) in $M_\R$ containing $0$ as an interior point. We define a notion of convexity with respect to $P$. 

\bd \label{d:adapt}
A homogeneous convex function $\varphi$ on $M_\R$ is {\em convex with respect to $P$}, if for each face $F$ of $P$ of positive dimension, the restriction $\varphi|_F$ has a unique minimum point in the interior of $F$. 
A   {\em \Kahler potential adapted to $P$} is a homogeneous degree two convex function $\varphi: M_\R \to \R$ that is convex with respect to $P$. 
\ed

\brem
A homogeneous convex function $\varphi$ on $M_\R$ is convex with respect to $P$, if the increasing sequence of level sets $\{ \varphi(\rho) < c\}$ meet the faces of $P$ in the interior first. 
\erem

\bp
For any convex polytope $P$ in $M_\R$ containing $0$ as an interior point, there exists a non-empty contractible set of \Kahler potential adapted to $P$. 
\ep
\bpf
First, we  show the existence of such potential $\varphi$. We will build the level set $S=\{\varphi=1\}$, and show that as we rescale $S$ to $\lambda S$, for $\lambda$ from $0$ to $\infty$, $S$ will meet the interior of each face $F$ first. We will proceed by first build a polyhedral approximation of $S$, then smooth it.

For each face $F$ of $P$, we pick a point $x_F$ in the interior of $F$ if $\dim F > 0$, or $x_F = F$ if $F$ is a point. Let $T$ be the simplicial triangulation of $P$ with vertices of $F$, then $T$ is also a barycentric subdivision of $P$. 
Let $\phi_T: P \to \R$ a piecewise linear convex function on $P$, with maximal convex domain the top-dimensional simplices of $T$, and such that for any $0 \leq d \leq n-1$, and any face $x_F$ of dimension $d$, $\phi_T(x_F) = c_d$ are the same for all such $F$. Such $\phi_T$ can be constructed inductively from $x_F$ with $\dim F$ from $0$ to $n-1$. Let $\phi_T$ be extended to $M_\R$ by linearity. Thus $\phi_T$ has a unique minium point in each face $F$. 

Let $\eta \in C^\infty_c(\R^n)$ be a bump function with $\int \eta = 1$, and let $\eta_\epsilon(x) = \eta(x/\epsilon)/\epsilon^n$. Let $\phi_{T,\epsilon} = \eta_\epsilon \star \phi_T + \epsilon |x|^2$, where $|x|$ is taken with respect to any fixed inner product on $\R^n$,   then $\phi_{T,\epsilon}$ is a linear combination of convex function hence still convex. Since $\phi_{T,\epsilon} \to \phi_T$ as $\epsilon \to 0$, for $\epsilon$ small enough, $\phi_{T,\epsilon}$ still  has a unique minimum point in each face $F$. And $S_{T, \epsilon} = \{\phi_{T,\epsilon}=1\}$ is a convex smooth boundary, such that $S_{T, \epsilon} \to S_T = \{\phi_T=1\}$ as $\epsilon \to 0$. 
Then, for small enough $\epsilon$, we can use $S_{T,\epsilon}$ as the contour of the homogeneous degree two convex function $\{\varphi(x)=1\}$. 

(2) Let $\kcal$ be the set of homogeneous degree two potential adapted to $P$.  Then there is surjective continuous map $\pi: \kcal \to \prod_{F , \dim F>0} \Int (F)$, by sending $\varphi$ to its critical points on each face. Since if two convex functions $\varphi_1, \varphi_2$ have the same critical points, then their convex linear combination $t \varphi_1 + (1-t)\varphi_2$ for $t \in [0,1]$ are still homogeneous degree two and with the same critical points, we see the fiber of map $\pi$ is convex hence contractible. Since the base of the fibration $Cr$ is contractible as well, we have $\kcal$ contractible. 

\epf

Let $P$ be a convex polytope in $M_\R$ containing $0$ as an interior point. Recall the definition of the dual polytope $P^\vee \subset N_\R$
\be P^\vee = \{ p \in N_\R \mid \la p, x \ra \leq 1 \; \forall x \in P \}. \ee
For any face $F \subset P$, there is dual face $F^\vee \subset P^\vee$, and $\dim_\R F + \dim_\R F^\vee = n-1$. 
We define three subsets of $M_\R \times N_\R$
\be L_{P} = \bigcup_F \cone(F) \times F^\vee , \quad L_{P^\vee} =\bigcup_F F \times \cone(F^\vee), \quad \Lambda_P  = \bigcup_F F \times F^\vee, \label{e:LP} \ee
where $F$ runs over the faces of $P$, and $\cone(F) = \R_{> 0} \cdot F$. 

\brem
 $L_P$ and $L_{P^\vee}$ are piecewise Lagrangians, and $\Lambda_P  = L_P \cap L_{P^\vee}$ is piecewise isotropic.  $L_P$ is the exterior conormal of $P^\vee$ in $T^*N_\R$, and $L^\vee_P$ is the exterior conormal of $P$ in $T^*M_\R$.  If we let $\varphi_{P, 1}$ be the piecewise linear function on $M_\R$, such that $P = \{ x: \varphi_{P, 1}(x) \leq 1 \}$, then $L_P$ morally is $\Gamma_{d \varphi_{P,1}}$. 
\erem 

\bl \label{l:dual-cone-ints}
Let $\varphi$ be a homogeneous degree two convex function on $M_\R$. $P, P^\vee$ be dual convex polytopes in $M_\R$ and $N_\R$ as above. Let $F$ be a face of $P$. Then there is a natural bijection 
\be \cone(F) \times F^\vee \cap \Gamma_{d \varphi}  \lra F \times \cone(F^\vee) \cap \Gamma_{d \varphi}. \ee 
\el
\bpf
If $(\lambda x, p) \in  \cone(F) \times F^\vee \cap \Gamma_{d \varphi}$, where $\lambda > 0$ and $x \in F, p \in F^\vee$, then by conic invariance of $\Gamma_{d \varphi}$, we have 
\be (x, p/\lambda) = \frac{1}{\lambda} (\lambda x, p) \in F \times \cone(F^\vee) \cap \Gamma_{d \varphi}. \ee
Sending $(\lambda x, p)$ to $(x, p/\lambda)$ is the desired bijection. 
\epf

Next, we give some equivalent characterization for convexity with respect to a polytope. 

\bp\label{p:equiv-conv}
Let $P$ be a convex polytope in $M_\R$ containing $0$ as an interior point. Let $\varphi$ be a homogeneous degree two convex function on $M_\R$.
The following conditions are equivalent: \\
(1) $\varphi$ is adapted to $P$. \\
(2) For each face $F$ of $P$, the smooth component $\Int(F \times \cone(F^\vee))$ of $L_{P}^\vee$  has a unique intersection with $\Gamma_{d \varphi}$. \\
(3) For each face $F$ of $P$, the smooth component $\Int(\cone(F) \times F^\vee)$ of $L_{P}$  has a unique intersection with $\Gamma_{d \varphi}$. 
\ep
\bpf
(2) is equivalent to (3) by Lemma \ref{l:dual-cone-ints}. 

(2) $\RA$ (1): since $\varphi|_F$ is still convex, hence as at most one minimum point in the interior, and any interior critical point is a minimum. Since 
\be \emptyset \neq F \times \cone(F^\vee) \cap \Gamma_{d \varphi} \subset T_F^* M_\R \cap  \Gamma_{d \varphi}  \ee
we see $\varphi|_F$ has a critical point.

(1) $\RA$ (2): for each face $F$ of $P$, let $x_F$ be the critical point of $\varphi|_F$, and let $H_F \subset M_\R$ be the affine hyperplane tangent to the contour of $\varphi$ at $x_F$. We claim that $H_F$ is a supporting hyperplane for $P$, and $P \cap H_F = F$. 
Then $p = d \varphi|_{x_F} \in T_{x_F}^* M_\R \cong N_\R$ is in the exterior conormal of $H_F$ (exterior with respect to $P$), hence  $p \in \cone(F^\vee)$. Thus, $(x, p) \in F \times \cone(F^\vee)$. 
\epf
%
%
%

A consequence of the proposition is the compatibility of the `adaptedness' with Legendre transformation. 
\bc \label{c:dual-adapt}
Let $P$ be a convex polytope in $M_\R$ containing $0$ as an interior point  and $P^\vee$ the dual polytope. Let $\varphi$ be homogeneous degree two convex function, and  $\psi$ the Legendrian dual of $\varphi$. Then $\varphi$ is adapted to $P$ if and only if $\psi$ is adapted to $P^\vee$. 
\ec

\section{Main Results \label{main}}
Let $Q$ be a convex lattice polytope in $N_\R$ containing $0= \alpha_0$. Let $\tcal$ be a coherent star triangulation of $Q$ based at $0$ with integral vertices, and $\pa \tcal$ be the subset of simplices not containing $0$. Let $A$ be the set of vertices of $\tcal$, and let $h: A \to \R$ induce $\tcal$ with $h(0)=0$. Fix a $\Theta: A \to T$ with $\Theta(\alpha_0)=\pi$. 

Let $\Pi_h \subset M_\R$ be the tropical amoeba of $(\tcal, h)$, and $P = C_{\alpha_0}$ be the connected component in $M_\R \RM \Pi_h$ corresponding to $\alpha_0$. 

Let $\varphi: M_\R \to \R$ be a homogeneous degree two convex function (i.e. $\varphi(\lambda x) = \lambda^2 \varphi(x)$ for all $\lambda>0$). \footnote{We may smooth $\varphi$ at an small neighborhood around $0 \in M_\R$, but this is irrelevant since we will use $\varphi$ only as $\varphi(\beta u)$ for $\beta \gg 1$ and $u$ in a neighborhood of $\pa P$. } We assume $\varphi$ is adapted to $P$, i.e, every face of $P$ contains a minimum of $\varphi$ in its interior. 

 Then the tropical localized polynomial is
\[ \wt f(z) := \sum_{\alpha \in A} e^{- \beta h(\alpha) - i \Theta(\alpha)}  \chi_{\alpha,\beta}(z) z^\alpha \] 
where   the monomial cut-off function $\chi_{\alpha, \beta}$ defined as in Definition  \ref{d:cutoff}. 
Denote the localized hypersurface and amoeba as $\wt \hcal := \wt f^{-1}(0)$ and $\wt \Pi:=\Log_\beta(\wt \hcal)$ respectively. We use $\wt C = \wt C_{\alpha_0}$ to denote the  complement labelled by $\alpha_0$. 

\bt\label{gradient}
The critical points for $\varphi|_{\pa \wt C}$ on the boundary of amoeba $\pa \wt C$ are indexed by simplices $\tau \in \pa \tcal$.  The critical point $\wt \rho_{\tau}$ for $\tau$  in $\pa \tcal$ has Morse index $\dim |\tau|$. The unstable manifold (downward flowing) $\wt W_{\tau}$ of $\wt \rho_{\tau}$ contains $\wt \rho_{\tau'}$ in its closure, if and only if $\wb \tau \supset \tau'$. 
\et

The Liouville structure of $\wt \hcal$ are induced from $( M_\CS, \omega, \lambda)$ (see Section \ref{ss:liouville}), where in coordinates
\[ \omega = \pa_{ij} \varphi(\rho) d \rho_i \wedge d\theta_j, \quad \lambda = \pa_i \varphi(\rho) d\theta_i. \]

The (candidate for) Liouville skeleton is defined as
\be \label{Liouville-skel}
\scal_{\beta, h, \Theta} = \bigcup_{\tau \in \tcal} (\beta \cdot \wt W_{\tau}) \times T_{\tau, \Theta} \subset M_\R \times M_T \cong M_\CS, \ee
where $\wt W_\tau$ is the unstable manifold from $\wt \rho_\tau$ and $T_{\tau, \Theta}$ is the subtorus of $M_T$ defined by 
\be \label{crit-torus}
 T_{\tau, \Theta} = \{\theta \in M_T : \la \theta, \alpha \ra = \Theta(\alpha), \text{for each vertex $\alpha$ in  $\tau$.} \} \ee

\bt\label{embed}
$\scal_{\beta, h, \Theta}$ is the Lagrangian Liouville skeleton for $(\wt \hcal, \omega|_{\wt \hcal}, \lambda|_{\wt \hcal})$
\et

The Lagrangian skeleton defined here can be related with the RSTZ skeleton via the 'projective' Legendre transformation $\Phi_\varphi^\infty: M_\R^\infty \isoto N_\R^\infty$, which is induced by homogeneous Legendre transformation $\Phi_\varphi: M_\R \to N_\R$. Let $q_M: M_\R \RM \{0\} \to M_\R^\infty$ and $q_N: N_\R \RM \{0\} \to N_\R^\infty$ be quotient by $\R_+$. Recall the RSTZ-skeleton $\Lambda^\infty_{\tcal, \Theta}$ is defined in the introduction \eqref{RSTZ-skeleton}. Let $id$ denote the identity map on $M_T$, then we have: 
\bt \label{compare}
\[ \Phi_\varphi^\infty \times id:  M_\R^\infty \times M_T \to N_\R^\infty \times M_T \]
induces a homeomorphism between $\scal_{\beta, h, \Theta}$ identified as  $(q_M \times id) (\scal_{\beta, h, \Theta})$ and $\Lambda^\infty_{\tcal, \Theta}$ identified as $(q_N\times id) (\Lambda^\infty_{\tcal, \Theta})$. 
\et

The main theorem then follows from Theorem \ref{embed} and \ref{compare}, and the diffeomorphism of $\hcal$ with $\wt \hcal$ from Proposition \ref{family}. 



\section{Gradient Flow: Proof of Theorem \ref{gradient}}
We will index the critical points by simplex $\tau$ in $\pa \tcal = \tcal \cap \pa Q$. We use $\tau_0$ to denote the simplex of $\conv(\{0\} \cup \tau)$. 
We will sometimes omit $\beta$ from the subscript to unclutter the notation. 

\ss{Convergence of smooth Convex Domain and Critical Points}
We fix an identification of $V \cong \R^n$ and take Euclidean metric on $V$ and the induced metric on $T^*V$ and $S^*V$.  We identify the sphere compactification boundary $T^\infty V = (T^*V - V)/ \R_{>0}$ with the unit cosphere bundle $S^*V$. If $U \subset V$ open  set with smooth boundary, then $S^*_UV$ is the one-sided unit conormal bundle of $\pa U$ with covectors pointing outward. The generalization to open convex set $U$ with piecewise smooth boundary is also straightforward. 

\bp \label{p:crit-conv}
Let $V \cong \R^n$ be a real vector space of dimension $n$, 
 $P \subset V$  a convex polytope containing the origin, $\varphi: V \to \R$ a potential adapted to $P$. 
Let $\{P_j\}$ be a sequence of convex bounded domains with smooth boundaries, such that the exterior conormals $L_j: = S^*_{P_j} V$ converges to $L:=S^*_P V$ in the cosphere bundle $ S^*V$ in the Hausdorff metric. Then for all large enough $j$, there is a one-to-one correspondence between face $F$ of $P$ and critical points of $\varphi$ on $\pa P_j$,  denoted as $\wt \rho_F$,  such that\\
(1) $\wt \rho_F$ has Morse index $n-1-\dim F$. \\
(2) As $\beta \to \infty$, $\wt \rho_F$ tends to the $\rho_F$, the minimum of $\varphi$ on the face $F$ of $P$. 
\ep

\bpf
(1) We express the critical point condition in terms of Legendrian intersection. Define the projection image of $\Gamma_{d \varphi}$ in $T^\infty V$ as
\be \label{Gdphi-1} \Gamma^\infty_{d\varphi} = (\R_{>0} \cdot \Gamma_{d\varphi}) / \R_{>0} \subset T^\infty V. \ee
Then $\Gamma^\infty_{d \varphi}$ is also the union of unit conormal for level sets of $\varphi$: 
\be \label{Gdphi-2} \Gamma^\infty_{d\varphi} = \bigcup_{c \in \R} S^*_{ \{\varphi(\rho)\leq c\}} V. \ee

The Legendrian $L = S^*_PV$ is a piecewise smooth $C^1$ manifold, where the smooth components $L_{F}$ are labelled by faces $F$ of $P$. If $\rho_F$ is a critical point of $\varphi$ on $F$, then there is a unique unit covector $p_F \in L_F$, such that $x_F = (\rho_F, p_F) \in L \pitchfork \Gamma^\infty_{d \varphi}$, and the intersection is transversal.

(2) Consider the unit speed geodesic flow $\Phi_R^t$ on the unit cosphere bundle $S^* V$.  Fix  any small flow time $1 \gg \epsilon>0$, since $\Phi^\epsilon_R: S^*V \to S^*V$ is a diffeomorphism, $\Phi^\epsilon_R(L_j)$ still converges to $\Phi^\epsilon_R(L)$ in Hausdorff metric. For any subset $A \subset V$, define 
\[ A^\epsilon := \{x: \dist(x,A)<\epsilon\} \] 
to be the $\epsilon$-fattening of $A$. If $A$ is a convex set, we have
$ \Phi_R^t(S^*_A V) = S^*_{A^\epsilon} V.$
Hence $\pa P^\epsilon$ is a $C^1$ hypersurface, and $\pa P_j^\epsilon \to \pa P^\epsilon$ in Hausdorff metric as $j \to \infty$. Define 
\[ L^t = \Phi_R^t(L), \quad  L^t_j = \Phi_R^t(L_j). \]

The geodesic flow applied to $\Gamma^\infty_{d \varphi}$ can be understood as follow
\[ \Phi_R^\epsilon(\Gamma^\infty_{d\varphi}) = \bigcup_{c \in \R} \Phi_R^\epsilon(S^*_{ \{\varphi(\rho)\leq c\}} V) = \bigcup_{c \in \R}  S^*_{ \{\varphi(\rho)\leq c\}^\epsilon} V.  \]
Define function $\wt \varphi^\epsilon$, such that $\{\wt \varphi^\epsilon(\rho) < c\} = \{\varphi(\rho)\leq c\}^\epsilon$, then $\wt  \varphi^\epsilon$ is a levelset convex function. By Lemma 2.7 of \cite{CE}, there exists a strictly increasing function $f: \R \to \R$, such that $\varphi^\epsilon = f \circ \wt \varphi^\epsilon$ is a convex function. 
Thus, we have 
\[  \Phi^\epsilon_R(\Gamma^\infty_{d \varphi}) = \Gamma^\infty_{d \varphi^\epsilon}  \quad \varphi^\epsilon \text{ is convex }. \]

Let $x_F^\epsilon = \Phi^\epsilon_R(x_F)$,  $\rho^\epsilon_F = \pi(x^\epsilon_F)$ in the expanded face $F^\epsilon = \pi(\Phi^\epsilon_R(L_F))$. Then $x^\epsilon_F$ is still the intersection of $\Gamma^\infty_{d \varphi^\epsilon}$ and $S^*_{P^\epsilon}V$, and $\rho_F^\epsilon$ is the unique Morse critical points of $\varphi^\epsilon$ restricted on $F^\epsilon$, and $\rho_F^\epsilon$ is in the interior of $F^\epsilon$. One may easily check that the Morse index of $\rho_F^\epsilon$ is $n-1-\dim F$. 

(3) We now prove that for large enough $j$, for each $F$, there is a unique critical points $\rho_{F,j}^\epsilon$ of $\varphi^\epsilon$ on $\pa P_j^\epsilon$ approaching $\rho_F^\epsilon$. 

Fix a small neighborhood $\wcal_F \subset \pa P^\epsilon$ near $\rho^\epsilon_F$, and for small enough $\delta$, let $\wt \wcal_F \cong \wcal_F \times (-\delta, \delta)$ be the flow-out of $\wcal_F$ under the Reeb flow for time in $(-\delta,\delta)$, with projection map $\pi_W: \wt \wcal_F \to \wcal_F$. We claim that for large enough $j$, $\pa P_j^\epsilon \cap \wt \wcal_F$ projects bijectively to $\wcal_F$, since otherwise this contradicts with $P^\epsilon_j$ being convex and the fiber of $\pi_W$ being straight-line segments Reeb trajectories. 
Thus, we have a sequence of smooth sections $\iota_j: \wcal_F \to \wt \wcal_F$ for large enough $j$, such that $\iota_j$ converges to the zero section in $C^1$.  

Let $f_j = \iota_j^* \varphi_\epsilon|_{\wt \wcal_F} \in C^\infty(\wcal_F, \R)$, and a smooth function $f_\infty = \iota_\infty^* \varphi_\epsilon|_{\wt \wcal_F}$, where $\iota_\infty: \wcal_F \into \wt \wcal_F$ is the identity map of zero section. Since $\iota_j \to \iota_\infty$ in $C^1$,  $f_j \to f_\infty$ in $C^1$. Since $f_\infty$ has a non-degenerate critical point, by stability of critical points under $C^1$-perturbation, $f_j$ has a unique critical point of the same index as $f_\infty$. 

(4) Finally, we show that there are no other critical points. Let $U_F$ be the preimage of $\wt W_F$ under $S^*V \to V$. 
Let $U$ be the union of all such $\wt U_F$. If $\delta>0$ is small enough, such that 
\[ \dist(\Gamma^{\infty}_{d\varphi^\epsilon} \RM U, L^\epsilon) > 3 \delta. \]
Then by the assumption that  $L_j^\epsilon$ converges to $L^\epsilon$ in Hausdorff metric,  we make take $j_0$ large enough, such that for all $j>j_0$ and all $x \in L_j^\epsilon$, $\dist(x, L^\epsilon) < \delta$. This shows 
\[ \dist(\Gamma^{\infty}_{d\varphi^\epsilon} \RM U, L_j^\epsilon) \geq \dist(\Gamma^{\infty}_{d\varphi^\epsilon} \RM U, L^\epsilon) - \dist(L_j^\epsilon, L^\epsilon) > 2 \delta, \]
hence there is no intersection between $L_j^\epsilon$ and $\Gamma^\infty_{d\varphi^\epsilon}$ away from $U$. 

(5) Since $\Phi_R^\epsilon$ is a diffeomorphism, the result about $L_j^\epsilon \cap \Gamma^\infty_{d\varphi^\epsilon}$ implies the same result about $L_j \cap \Gamma^\infty_{d\varphi}$, and we finish the proof of the proposition. 
\epf

\subsection{Proof of Theorem 1: Critical Points and Unstable Manifolds}
\bp\label{p:crit-beta}
For large enough $\beta$, there is a one-to-one correspondence between simplices $\tau \in \pa \tcal$ and critical points of $\varphi$ on $\pa \wt C$,  denoted as $\wt \rho_\tau$,  such that\\
(1) $\wt \rho_\tau$ has Morse index $\dim \tau$. \\
(2) As $\beta \to \infty$, $\wt \rho_\tau$ tends to the $\rho_\tau$, the minimum of $\varphi$ on the face $\tau^\vee_0$ of $P$. 
\ep
\bpf
First we approximate $\pa \wt  C $ by its convex model $\wh C$ (see Definition \ref{d:convexC}). Then by  Proposition \ref{p:crit-conv}, we have critical points $\{ \wh \rho_{\tau} \}$ on $\wh C$ indexed by $\tau \in \pa \tcal$. Then, a perturbation argument shows $\varphi$ has critical points on $\pa \wt  C $ as $\{ \wt \rho_{\tau} = \wh \rho_{\tau}\}$. 
\epf

Next, we prove that the unstable manifold $\wt W_{\tau}$ for critical point $\wt \rho_{\tau}$ are cells of a dual polyhedral decomposition of $\pa P$. This is true not only in the combinatorial sense, but in a more refined geometrical sense. 

\bp
For large enough $\beta$, and for any $\tau \in \pa \tcal$,  the unstable manifold $\wt W_{\tau}$ is a smooth manifold of dimension $\dim \tau$, and $\wt \rho_{\tau'} \in \pa{\wt W_{\tau}} $ if and only if $\tau' \subset \tau$.  
\ep
\bpf
The statement of $\dim \wt W_\tau$ follows from the Morse index of $\wt \rho_\tau$. For any critical point $\rho_{\tau}$, take a small enough ball $B$ of radius $\epsilon$ around it, then $B$ can be stratified by the limit of gradient flow. For each facet $\sigma^\vee$ of the polytope $P$ adjacent to $\tau^\vee$, there is an open ball $U_\sigma$ in $\pa B$ whose points flows to critical point $\rho_{\sigma}$. If a face $\tau'^\vee$ adjacent to $\tau^\vee$ can be written as $\tau'^\vee = \sigma^\vee_1 \cap \cdots \cap \sigma^\vee_k$ for facets $\sigma^\vee_i$, then points in the relative interior of $ \cap_{i=1}^k \wb{U_{\sigma_i}}$ will flow to $\rho_{\tau'}$. 
\epf

Now we give an explicit description of the unstable manifold. Let $\Phi^\infty_{\varphi}$, $q_N$ and $q_M$ be as in the statement of Theorem \ref{compare}.  
 
\bp \label{p:cone}
For all large enough $\beta$, and $\tau \in \pa \tcal$,  $\Phi_\varphi^\infty$ induces a homeomorphism 
\[ \Phi_\varphi^\infty: q_M(\wt W_\tau) \isoto q_N(\tau). \]
\ep
\bpf
Without loss of generality, we set $h(0)=0$. We will first do the computation on the convex model $\wh C$, then state the necessary modifications for $\wt C$. 

Then the defining function for $\wh C$ can is 
\[ 1 = \sum_{\alpha \in A \RM \{0\}} e^{\beta l_\alpha(u)} \chi(\beta l_\alpha(u) + \sqrt \beta) = :\wh F(u) \]
For any point $u \in \wh C$, we define the simplex
\[ \tau(u) = \conv \{\alpha \in A \RM \{0\} : \chi(\beta l_\alpha(u) + \sqrt \beta) > 0\}. \]

Then the gradient of $\varphi$ on $\wh C$ can be expressed as 
\[ \nabla(\varphi|_{\wh C}) = \nabla \varphi - c_1 \nabla \wh F \]
where $c_1 = \frac{\la \nabla \varphi, d \wh F \ra}{\la \nabla \wh F, d \wh F \ra }$. Since by Proposition \ref{p:gradphi}, $\nabla \varphi$ is in the outward radial direction, and $\wh F$ is a convex function with bounded sub-level set, hence $\la \nabla \varphi, d \wh F \ra > 0$. Combining $\la \nabla \wh F, d \wh F \ra > 0$, we have $c_1 > 0$. 

For $u$ on the unstable manifold $\wh W_\tau$, we have $\tau(u) \subset \tau$. The defining function $\wh F$ for a neighborhood of $u$ can be written as 
\[ \wh F_\tau(u) = \sum_{\alpha \in \tau} e^{\beta l_\alpha(u)}  \chi(\beta l_\alpha(u) + \sqrt \beta) \]
thus
\[ d \wh F_\tau = \sum_{\alpha \in \tau} (e^{x} \chi(x+\sqrt{\beta}))'|_{x = \beta l_\alpha(u)} \cdot \alpha \in \Int \cone(\tau). \]
And at the critical point 
\[ d \varphi(\wh \rho_\tau) = c_1 d \wh  F_\tau \in \Int \cone(\tau). \]

If $\gamma: (-\infty, +\infty) \to \wh C$ is an integral curve for $-\nabla(\varphi|_{\wt C})$ with $\lim_{t\to -\infty}\gamma(t) = \wh \rho_\tau$, then under Legendre transformation we have a curve $ \eta(t)$, such that 
\[ \lim_{t \to -\infty} \eta (t) = d \varphi(\wt \rho_\tau) \in \Int \cone(\tau)\]
and using Lemma \ref{p:gradf} and Proposition \ref{p:gradphi}
\[ \frac{d}{dt} \eta(t) = (\Phi_\varphi)_*(-\nabla(\varphi|_{\wt C})) = (\Phi_\varphi)_*(-\nabla \varphi + c_1 \nabla F) \in \R (p \pa_p) +   \Int \cone(\tau), \]
where $p \pa_p$ is the radial vector field on $N_\R$. Thus $\eta(t)$ is within the cone $\Int \cone(\tau)$ for all $t \in \R$. This shows that 
\[ \Phi_\varphi^\infty(q_M(\wh W_\tau)) \subset q_N(\tau). \]
Using induction on dimension of $\tau$ from $0$ to $n-1$, we can show the image is onto. 

Now consider $\wt C$. One need to replace $\chi(\beta l_\alpha(u) + \sqrt \beta)$ by $\chi_{\alpha, 0})(u)$ in defining $\wt F$. And one can still show that $d \wt F_\tau(u) \in \Int \cone(\tau)$, the rest is the same as $\wh C$. 
\epf

\section{Liouville Flow: Proof of Theorem \ref{embed} and \ref{compare}}
First we find all the critical points (manifolds) of the Liouville vector field on $\wt \hcal$. We show that they are exactly the preimage of critical points of $\varphi |_{\pa \wt C}$ under $\Log_\beta$, which are tori of various dimensions. The more difficult part is to show there are no other critical points. 

Then, we study the Liouville flow trajectory from these critical manifolds. There are two key points:
\begin{enumerate}
\item We write the fiberwise Liouville vector field as the ambient Liouville vector field in $M_\CS$ subtract its symplectic orthogonal component, then show that on the 'positive loci' (Definition \ref{posi}), the symplectic orthogonal component is  proportional to the Hamiltonian vector field $X_{\Im \wt f}$. 
\item We show that the unstable manifold correponding to the critical manifold indexed by $\tau \in \pa \tcal$, is geometrically identified with the simplex $\tau$ under the projective Legendre transformation $\Phi_\varphi^\infty: M^\infty_\R \to N^\infty_\R$. This determines the unstable manifolds. 
\end{enumerate}

Recall our notations: 
\begin{itemize}
\item The complex hypersurface $\hcal$ defined by $f(z)=0$: 
\[ f(z) = \sum_{\alpha \in A} f_\alpha(z) = \sum_{\alpha \in A} z^\alpha e^{-\beta h(\alpha) - i \Theta(\alpha)}.  \]
\item The tropical localized hypersurface $\wt \hcal$, defined by $\wt f(z)=0$: 
\[ \wt f (z) =\sum_{\alpha \in A} \wt f_\alpha(z) = \sum_{\alpha \in A} f_\alpha(z) \chi_\alpha(u), \]
\item The real-valued functions,  $F_\alpha(z) = |f_\alpha(z)|$ and 
\[ \wt F(z) = -1 +  \sum_{0 \neq \alpha \in A}  F_\alpha(u) \chi_\alpha(u). \]
\end{itemize}

The proof of Theorem 2 follows from Propositions \ref{unstable}, \ref{lag} and \ref{no-other}. The proof of Theorem 3 follows from Propositions \ref{p:cone} and \ref{unstable}. 

\subsection{Liouville Vector Field}
Take any point $z \in \wt \hcal$, we have 
\[ X_\lambda(z) = X_\lambda^\|(z) + X_\lambda^\perp(z). \]
where $X_\lambda^\perp(z)$ is symplectically orthogonal to $T_z\wt \hcal$. 
We note that $X_\lambda^\|(z) = X_{\lambda_\hcal}(z)$, since for any $v \in T_z \wt \hcal$, 
\[ \omega_\hcal(X_\lambda^\|(z), v) = \omega(X_\lambda(z) - X_\lambda^\perp(z), v) = \omega(X_\lambda(z) - X_\lambda^\perp(z), v)  = \lambda(v) = \lambda_\hcal(v).\]
And $X_\lambda^\perp(z)$ is the symplectic horizontal lift of $\wt f_*(X_\lambda(z)) \in T_0 \C$. 

\bd\label{posi}
The {\em positive loci} $\wt \hcal^+$ is the subset of $\wt \hcal$  where $\wt f_0 = -1$ and $\wt f_\alpha \geq 0$ for all  $\alpha \neq 0$.  
\ed

\brem
An equivalent definition is that 
$ \wt \hcal^+ = \Log_\beta^{-1}(\pa \wt C)$. 
\erem

\bp
For all  $z \in \wt \hcal^+$, we have $X_{\Im \wt f}(z)$ positively proportional to $X_\lambda^\perp(z)$.  
\ep
\bpf
Since both vectors are symplectic orthogonal to $\wt \hcal$, we only need to check their image under $\wt f_*$ are positively proportional to each other. 

First we study   $X_{\Im \wt f}(z)$. On $\wt \hcal^+$, we have $\wt f = \wt F$. We also have 
\bea
d \wt f &=&  \sum_{\alpha \in A} \wt F_\alpha(u) \la \alpha, d(\rho+i\theta) \ra + \sum_{\alpha \in A} f_\alpha(z) d \chi_\alpha(u) \\
&=& d \wt F(u) + i \sum_{\alpha \in A} \wt F_\alpha(z) \la \alpha, d\theta \ra
\eea
Hence 
\[ d (\Im f) = \Im df = \sum_{\alpha \in A} \wt F_\alpha(u) \la \alpha, d\theta \ra \]
Thus 
\be \label{ximf} X_{\Im f} = \sum_{\alpha \in \pa A} \wt F_\alpha(u) \sum_{i,j}\alpha_i g^{ij}(\rho) \pa_{\rho_j} \ee
compare with 
\[\nabla (\wt F) =  g^{-1}(d \wt F) = g^{-1}( \sum_{\alpha \in \pa A} F_\alpha(u) \chi_{\alpha}(u) \la \alpha, d\rho\ra + F_\alpha d\chi_\alpha) = X_{\Im f}+ O(e^{-\sqrt{\beta}}). 
 \]
We thus have
\[ \la d\wt f, X_{\Im \wt f} \ra = \la d \wt F, X_{\Im \wt f} \ra  = \|\nabla \wt F\|^2 + O(e^{-\sqrt{\beta}}) > 0. \]

Next, we study $X_\lambda^\perp(z)$. We have 
\[ \la d \wt f, X_\lambda^\perp(z) \ra =  \la d \wt f, X_\lambda (z) \ra = \la d \wt f, \nabla \varphi \ra = \la d \wt F, \nabla \varphi \ra \]
Since $\nabla \varphi$ is positively proportional to the radial vector field $u\pa_u$ by Proposition \ref{p:gradphi}, and $\la d\wt F,  u\pa_u\ra > 0$. We have also $\la d \wt f, X_\lambda^\perp(z) \ra>0$. 

Since $\wt f_*(X_\lambda^\perp(z))$ and $\wt f_*(X_{\Im f})$ are both in the positive direction of $T_0 \C$, $X_\lambda^\perp(z)$ is positively proportional to $X_{\Im \wt f}(z)$.
\epf

\subsection{Critical Manifolds}
Recall  from the previous section, that on the boundary of the amoeba $\pa \wt C$, the critical points of $\varphi$ are indexed by $\tau \in \pa \tcal$ as $\wt \rho_\tau$.

\bp
The preimages $\Crit_\tau:= \Log_\beta|_{\wt \hcal}^{-1}(\wt \rho_\tau)$ are critical manifolds. 
\ep
\bpf
Since the critical points $\wt \rho_\tau$ are in the 'good' region $\wt \hcal^{good} \subset \wt \hcal$, where the monomial cut-off functions $\chi_{\alpha}$ are either zero or one, hence the hypersurface $\wt \hcal^{good}$ is holomorphic. Thus, zero of $d\varphi|_{\wt\hcal}$ is also zero of $d^c\varphi|_{\wt \hcal}$. 
\epf

\bp
For each $\tau \in \pa \tcal$, $\Log_\beta|_{\wt \hcal}^{-1}(\wt \rho_\tau) = \{\beta \wt \rho_\tau\} \times T_{\tau, \Theta}$, where  $T_{\tau, \Theta}$ is defined in \eqref{crit-torus}. 
\ep
\bpf
Since $\wt \rho_\tau$ is on the boundary $\pa \wt C$, we have 
$ 1 = \sum_{\alpha \in \tau } F_\alpha(z)$. 
Comparing with the defining equation of $\wt \hcal$ in a neighborhood of $\Log_\beta|_{\wt \hcal}^{-1}(\wt \rho_\tau)$, we have
$ 1 = \sum_{ \alpha \in \tau} f_\alpha(z).$
Hence 
$ 0 = \arg(f_\alpha(z)) = \la \alpha, \theta\ra - \Theta(\alpha)$
 for each vertex $\alpha$ in $\tau$. Thus the fiber is the torus $T_{\tau, \Theta}$. 
\epf

\subsection{Unstable Manifolds}

\bp \label{notheta}
The Liouville vector field $X_{\lambda_\hcal}$ on the positive loci $\wt \hcal^+$ does not change the $\theta$ coordinate.  
In particular, the positive loci $\wt \hcal^+$ is preserved under the Liouville flow. 
\ep
\bpf
Since
$ X_{\lambda_\hcal} = X_\lambda ^\| = X_\lambda - X_\lambda^\perp$, 
suffice to check that $X_\lambda$ and $ X_\lambda^\perp$ does not change $\theta$ coordinates. We have 
$X_\lambda \propto \rho \pa_\rho$, and $X_\lambda^\perp \propto X_{\Im \wt f}$.  From Eq. \eqref{ximf}, we see $X_{\Im \wt f}$ has no $\theta$-component. Hence $\la X_{\lambda_\hcal}, d\theta \ra = 0$. 
\epf

\bp \label{unstable}
For any $\tau \in \pa \tcal$, the unstable manifold for $\Crit_\tau$ is $\wt W_{\tau} \times T_{\tau, \Theta}$. 
\ep
\bpf
From Proposition \ref{notheta}, we see the flowout of $\Crit_\tau$ by the Liouville flow does not affect the $M_T$ component. Thus Louville flow $X_{\lambda, \wt \hcal}$ on $\wt \hcal$ induces a flow on $\wt \hcal^+$, and it descend to $\pa \wt C$, for which we  denote as $X_{\lambda, \pa \wt C}$. 

On $\pa \wt C^{good}$, $X_{\lambda, \pa \wt C}$ agrees with $\nabla (\varphi|_{\pa \wt C})$. And they have the same critical points set. On $\pa \wt C^{bad}$, we have 
\[ \|X_{\lambda, \pa \wt C} - \nabla (\varphi|_{\pa \wt C}) \| = O(e^{-\sqrt{\beta}}).  \]
Despite individual flowlines for the two vector fields with the same starting point in the good region may be split after flow through a bad region, we claim that for each critical point $\wt \rho_\tau$, the unstable manifolds $\wt  W^{X_\lambda}_{\tau}$ and $\wt W^{\nabla \varphi}_{\tau}$ for the two flows  are the same. 

Let $\tau \in \pa \tcal$ have vertices $\{\alpha_1, \cdots, \alpha_k\}$. Then 
\[ \nabla(\varphi|_{\wt C}) = \nabla \varphi - c_1 \nabla \wt F \in \R \cdot \rho \pa \rho + (\Phi_\varphi)_*^{-1} (\Int \cone \tau) \]
and 
\[ X_{\lambda, \pa \wt C}= X_\lambda - X_\lambda^\perp = X_\lambda - c(u) X_{\Im \wt f} \in \R \cdot \rho \pa \rho + (\Phi_\varphi)_*^{-1} (\Int \cone \tau) \]
where we used $X_\lambda^\perp $ positively proportional to $X_{\Im \wt f}$, and $X_{\Im \wt f}$ is given by Eq. \eqref{ximf}. By similar argument in Proposition \ref{p:cone} that $W^{\nabla \varphi}_{\tau}$ is dual to $\tau$ via $\Phi_\varphi^\infty$, we have $W^{X_\lambda}_{\tau}$ is dual to $\tau$ via $\Phi_\varphi^\infty$. Thus $\wt W^{\nabla \varphi}_{\tau}$ and $\wt W^{\nabla \varphi}_{\tau}$ has to be the same. We drop the superscripts and denote both as $\wt W_\tau$. 
\epf

\bp \label{lag}
For each $\tau \in \pa \tcal$, the unstable manifold $\wt W_\tau \times T_{\tau, \Theta}$ is a Lagrangian in $\wt \hcal$.
\ep
\bpf
One can use the property of the Liouville flow to show the unstable manifold is isotropic, and then counting dimension
\[ \dim_\R \wt W_\tau \times T_{\tau, \Theta} = (\dim \tau) + n-(\dim \tau+1) = n-1 = \half \dim_\R \wt \hcal. \] 
We give an alternative proof. By Proposition \ref{p:cone}, we have 
\[ \Phi_{\varphi} (\cone \wt W_\tau) \times T_{\tau, \Theta}= \cone \tau \times  T_{\tau, \Theta}.\]
However, $\cone \tau \times  T_{\tau, \Theta}$ is part of the conormal Lagrangian $T^*_{T_{\tau, \Theta}} M_T$ for the submanifold $T_{\tau, \Theta}$ in $M_T$. Since $\Phi_{\varphi} \times id$ is a symplectomorphism between $M_\CS$ and $T^* M_T$, we get $\cone \wt W_\tau \times T_{\tau, \Theta}$ is a conical Lagrangian in $M_\CS^*$. Finally, a Lagrangian restricts to a symplectic submanifold is isotropic. Thus by dimension counting,
\[ \wt W_\tau \times T_{\tau, \Theta} = \left( \cone \wt W_\tau \times T_{\tau, \Theta} \right) \bigcap \wt \hcal \]
is a Lagrangian in $\wt \hcal$. 
\epf

\subsection{No other Critical Points}

\bp \label{no-other}
There are no other zero of the Liouville vector field away from $\{\Crit_\tau\}$. 
\ep
\bpf
Suffice to prove that there are no zero of the Liouville vector field outside of the positive loci $\wt \hcal^+$. 
Here we only give the sketch the proof. We look at the good region first. Then $d^c \varphi|_\hcal = 0$ is equivalent to $d\varphi|_\hcal=0$, we only need to check there are no critical point for $\varphi$. 

Suppose there is a critical point $\varphi$ at $z \in \wt \hcal^{good}$, the terms labeled by $\alpha_1, \cdots, \alpha_k$ are non-zero, i.e., near $z$, $\wt \hcal$ is defined by 
\[ \sum_{i=1}^k f_{\alpha_i}(z) = 0. \]
Let $\tau \in \tcal$ be the simplex with vertices $\{\alpha_1, \cdots, \alpha_k\}$. Let $\tau^\vee$ be the cell in the tropical amoeba $\Pi$, and $U_\tau \subset \wt \hcal^{good}$ where the defining equation is as above. We split in to two cases below. Recall $P$ is the  polytope corresponding to vertex $0 \in \tcal$. Let $g_0$ denote the Euclidean metric on $M_\R$ after identification $M_\R \cong \R^n$. 

{\bf (1) The case $0 \notin \tau$.} Then $\tau^\vee$ is a non-compact cell in $\Pi$, and intersect the amoeba polytope $P$ at face $F_\tau = P \cap \tau^\vee$.  

Let $u = \Log_\beta(z)$, and let $u'$ denote the orthogonal projection w.r.t $g_0$ to the cell $\tau^\vee_0$. Then $dist_{g_0} (u, u') = O(1/\sqrt{\beta})$. Let $u''$ denote the minimum of $\varphi$ on $F_\tau$. We claim that 
\[ \varphi(u'') < \varphi(u'), \]
since the increase level set of $\varphi$ meet the convex cell $\tau^\vee$ first at $u''$. 

Let $v = u'' - u' \in M_\R$. If we view $v$ as a tangent vector at $u''$, then $\la d\varphi(u'), v \ra < 0$. Since $u$ and $u'$ are $O(1/\sqrt{\beta})$ close, we also have $\la d\varphi(u'), v \ra < 0$. Finally, one can check $v$ can be lifted as a tangent vector to $T_z \wt \hcal$, hence $d\varphi \neq 0$ at $z$. 

{\bf (2) The case $0 \in \tau$.} Without loss of generality, we may assume $\Theta(0)=\pi, h(0)=0$, and $\alpha_{k}=0$. Thus, the defining equation of $\wt \hcal$ near $z$ can be written as 
\[ 1 = \sum_{i=1}^{k-1} f_{\alpha_i} = \sum_{i=1}^{k-1} e^{-i \Theta(\alpha_i) - \beta h(\alpha_i)} e^{\beta \la \alpha_i, u\ra + i \la \alpha_i, \theta \ra } =: F(u, \theta) \]
Suppose $z$ is a critical point of $\varphi|_{\{F=1\}}$, then there exists $c_1, c_2$, such that 
\[ d \varphi(\rho) = c_1 d \Re F(\rho, \theta) + c_2 d \Im F(\rho, \theta). \]
However, since $d \varphi(\rho)$ has no $d \theta$ component, hence $d\theta$ on the RHS need to be cancelled out. Using all the $\alpha_i$ are linearly independent, we can check this is only possible if all $\arg(f_{\alpha_i})$ are equal or differ by $\pi$. Since $\sum_i f_{\alpha_i}=1$, we get all $f_{\alpha_i} \in \R$, and at least one is positive. 

If all of $f_{\alpha_i}(z)$ are positive, then there is nothing to show, since we want to prove all the critical points lies on the positive loci. 

If not all of $f_{\alpha_i}(z)$ are positive, say for $i=1,\cdots,m$, $f_{\alpha_i}(z)<0$, then $u$ lies on the real hypersurface
\[ 1 = - e^{\beta l_{\alpha_1}(u)} - \cdots -  e^{\beta l_{\alpha_m}(u)} + \cdots + e^{\beta l_{\alpha_{k+1}}(u)} =: H(u).  \]
near the face $\tau^\vee$ on $P$. If we further require $d\varphi$ to be in the $\R$-span of $\alpha_1, \cdots, \alpha_{k-1}$, then $u$ has to be near the critical point of $\varphi$ on face $\tau^\vee$. One can show that $d\varphi$ has to be in the $\R_+$-span of $\alpha_1, \cdots, \alpha_{k-1}$. Hence, there does not exists $c \in \R$, such that $d \varphi(u) = c d H(u)$.

This concludes the discussion for $z$ in the good region. If $z$ is in the bad region, where at least one $0 < \chi_{\alpha}(z) < 1$, we will approximate the bad region using good region in the following way. Define a different set of cut-off functions, by changing the cut-off threshold from $-\sqrt{\beta}$ to $-2\sqrt{\beta}$, ie. redefine 
 \[  \chi_{\alpha, \alpha', \beta}(u)=\chi(\beta(l_{\alpha}(u)  - l_{\alpha'}(u)) + 10 \sqrt{\beta}) \]
 in Definition \ref{d:cutoff}. Denote the new tropical localized hypersurface $\wt \hcal_{10}$. We claim the Hausdorff distance between $\wt \hcal$ and $\wt \hcal_{10}$ in $M_\CS$ is $O(e^{-c\sqrt{\beta}})$ for some $c>0$. Furthermore, their unit conormal bundles $S^*_{\wt \hcal} M_\CS$ and $S^*_{\wt \hcal_{10}} M_\CS$ should have distance $O(e^{-c\sqrt{\beta}})$ as well. A zero of $d^c (\varphi|_{\wt \hcal})$ corresponds to an intersection of $\Gamma_{d^c \varphi}^\infty \subset S^* M_\CS$ with $S^*_{\wt \hcal} M_\CS$, where
 \[ \Gamma_{d^c \varphi}^\infty = (\Gamma_{d^c\varphi} \cap \dot T^*(M_\CS) ) / \R_+ \subset T^\infty (M_\CS) \cong S^* (M_\CS). \]  
(cf Definition \eqref{Gdphi-1} and \eqref{Gdphi-2}). Thus the bad region of $\wt \hcal$ can be approximately by part of good regions in $\wt \hcal_{10}$, where we know there does not exists critical points of $\varphi$ away from the positive loci, hence there are no critical points in the bad region of $\wt \hcal$ away from the positive loci. 
\epf

\end{document}